\documentclass[a4paper,12pt]{article}
\usepackage{amsthm}
\usepackage{amsmath,amsfonts,amssymb}
\usepackage{a4wide}
\usepackage[usenames,dvipsnames]{color}
\usepackage[all,cmtip]{xy}
\usepackage[verbose,colorlinks=true,linktocpage=true,linkcolor=blue,citecolor=blue]{hyperref}
\usepackage{footnote}
\usepackage{tikz}
\usepackage{extarrows}
\usepackage[title]{appendix}
\usepackage{marvosym}
\usepackage{ulem}

\usepackage{indentfirst} 
\usepackage{geometry}
\newtheorem{theo}{{Theorem}}[section]
\newtheorem{lemm}[theo]{Lemma}
\newtheorem{rema}[theo]{Remark}
\newtheorem{defi}[theo]{Definition}
\newtheorem{coro}[theo]{Corollary}
\newtheorem{prop}[theo]{Proposition}

\numberwithin{equation}{section}

\hypersetup{
    colorlinks=true, 
    linkcolor=blue, 
    urlcolor=blue, 
}

\allowdisplaybreaks[1]

\begin{document}

\title{The double super Yangians in type A for arbitrary $0^m 1^n$-sequences and their bosonic representations}

\author{ Pengfa Xu${}^1$, Hongda Lin${}^2$ and Honglian Zhang${}^{1,}$\thanks{Corresponding Author.~~Email:~hlzhangmath@shu.edu.cn}}

\date{}

\maketitle

\begin{center}
\footnotesize
\begin{itemize}
\item[1] Department of Mathematics, Shanghai University, Shanghai 200444, PR~China.
\item[2] Shenzhen International Center for Mathematics, Southern University of Science and Technology, Shenzhen 518055, PR~China.

\end{itemize}
\end{center}
\begin{abstract}
In this paper, we introduce the double super Yangian
$\mathrm{DY}_{h}(\mathfrak{gl}_{m|n}^{\mathfrak{s}})$ and $\mathrm{DY}_{h}(\mathfrak{sl}_{m|n}^{\mathfrak{s}})$ associated with any fixed $0^{m}1^{n}$-sequence $\mathfrak{s}$. First, we establish an explicit isomorphism between the Drinfeld and R-matrix presentations of $\mathrm{DY}_{h}(\mathfrak{gl}^{\mathfrak{s}}_{m|n})$. We then generalize the notion of the quantum Berezinian to $\mathrm{DY}_{h}(\mathfrak{gl}_{m|n}^{\mathfrak{s}})$, and employ it to construct the R-matrix presentation of $\mathrm{DY}_{h}(\mathfrak{sl}_{m|n}^{\mathfrak{s}})$ and prove that it is isomorphic to the Drinfeld presentation. As an application, we present level-1 bosonic representations for  $\mathrm{DY}_{h}(\mathfrak{gl}_{m|n}^{\mathfrak{s}})$ and $\mathrm{DY}_{h}(\mathfrak{sl}_{m|n}^{\mathfrak{s}})$ in terms of their Drinfeld current generators.

\end{abstract}

\medskip
\section{Introduction}
Quantum groups, introduced by Drinfeld \cite{D87} and Jimbo \cite{Ji85}, have been extensively studied over the past few decades and now play important roles across mathematics and physics. A key example is the Yangian $\mathrm{Y}(\mathfrak{g})$, a canonical deformation of the universal enveloping algebra $\mathrm{U}(\mathfrak{g}[t])$ for a finite-dimensional semisimple Lie algebra $\mathfrak{g}$. The general linear Yangian $\mathrm{Y}(\mathfrak{gl}_{n})$ was first given by Tarasov \cite{Ta84} for $\mathfrak{gl}_{2}$ and later extended to the general case by Drinfeld \cite{D85}. It is a unital associative algebra that also admits a presentation in terms of defining relations written in a specific matrix form within the R–matrix formalism developed by Faddeev, Reshetikhin, and Takhtajan \cite{FRT90}. This alternative presentation endows the Yangian with a natural bialgebra structure, which can be extended to a Hopf algebra by adjoining an antipode.

In 1988, Drinfeld \cite{D88} proposed a new current realization for the Yangian, known as the Drinfeld new presentation, and further conjectured its equivalence to the R-matrix presentation without proof.
The conjecture was first confirmed by Brundan and Kleshchev \cite{BK05} for type $\boldsymbol{A}$, using the Gauss decomposition of the generator matrix $T(u)$. In 2018, the equivalence of the two presentations for types $\boldsymbol{B},\boldsymbol{C},\boldsymbol{D}$ was established in \cite{JLM18}.

In the super case, Nazarov \cite{N91} defined the super Yangian associated with the Lie superalgebra $\mathfrak{gl}_{m|n}$ as a super analogue of $\mathrm{Y}(\mathfrak{gl}_{n})$ via the R-matrix presentation. Subsequently, Gow \cite{G07} obtained a Drinfeld-type presentation corresponding to a standard $01$-sequence, and Peng \cite{P16} constructed parabolic presentations of the super Yangian $\mathrm{Y}(\mathfrak{gl}_{m|n}^{\mathfrak{s}})$ for an arbitrary $01$-sequence $\mathfrak{s}$. A further development came from Tsymbaliuk \cite{T20}, who linked arbitrary Drinfeld super Yangians of type $\boldsymbol{A}$ via odd reflections, which correspond to index permutations. Finite-dimensional irreducible representations of $\mathrm{Y}(\mathfrak{gl}_{m|n}^{\mathfrak{s}})$ were  described via these odd reflections \cite{L22, M22} in 2022. In addition, the Drinfeld presentation for orthosymplectic Yangians was given in \cite{M24}, building upon the R-matrix presentation introduced in \cite{AACFR03}.

The quantum double, obtained by pairing a quantum group with its Hopf dual \cite{D88}, provides the conceptual basis for the universal R-matrix of quantum enveloping algebras \cite{CP94}. This construction solves the quantum Yang-Baxter equation on every irreducible representation and drives advances in representation theory \cite{BFR92, EFM01}, $q$-conformal field theory \cite{EFR96}, and knot theory \cite{R90, RT91}. The double Yangian, defined as the quantum double of the Yangian $\mathrm{Y}(\mathfrak{g})$, applies naturally to massive field theory. The case without central extension was first considered for $\mathfrak{g}=\mathfrak{sl}_{2}$ \cite{K97}. The general case for $\mathfrak{g}=\mathfrak{gl}_{n} , \mathfrak{sl}_{n}$ was subsequently investigated by Iohara \cite{I96} within the R-matrix formalism with central extensions, where he also derived the Drinfeld commutation relations and constructed representations. Recently, Jing and Yang \cite{YJ24} established the isomorphism between the R-matrix and Drinfeld presentations of the double Yangian in type $\boldsymbol{A}$ and constructed the central elements at the critical level, while the analogous isomorphism for types $\boldsymbol{B}, \boldsymbol{C}, \boldsymbol{D}$ was given in \cite{JYL20}.

The super counterpart of the double Yangian was also formulated via the R-–matrix formalism, first introduced by Zhang \cite{Z97} for the standard general linear super case. Very recently, the quantum Berezinian, the PBW theorem, and reflection algebras for the double Yangian of $\mathfrak{gl}_{m|n}$ associated with the standard $01$-sequence were established in \cite{BK25-1, BK25-2}.

Vertex representations of quantum affine algebras of simply laced types were constructed by Frenkel and Jing in \cite{FJ88}, while the corresponding vertex operator representations and Drinfeld realizations for simply laced quantum affine superalgebras have  appeared in \cite{XZ}. For double Yangians, level 1 bosonization and vertex operators for type $\boldsymbol{A}$ were constructed in \cite{I96}, higher-level bosonization of $\mathrm{DY}_{h}(\mathfrak{sl}_{n})$ was obtained in \cite{DHHZ98}, and level-1 bosonization for types $\boldsymbol{B},\boldsymbol{C},\boldsymbol{D}$ was established in \cite{JYL20}. Using the results of \cite{I96}, Kozic constructed commutative operators for the double Yangian $\mathrm{DY}(\mathfrak{sl}_{n})$ in \cite{Ko18}.

In type $\boldsymbol{A}$, 01-sequences are commonly used to specify the parities of generators of (affine) Lie superalgebras, and distinct parity sequences correspond to different Dynkin diagrams. A distinctive feature of Lie superalgebras is that they admit multiple non-isomorphic Dynkin diagrams. The isomorphism problem for Lie superalgebras of finite or affine type associated with distinct Dynkin diagrams was resolved in \cite{LSS}. Subsequently, the analogous questions for quantum finite/affine superalgebras and type $\boldsymbol{A}$ super Yangians were addressed by Yamane \cite{Y99} and Tsymbaliuk \cite{T20}, respectively. Nevertheless, the corresponding problem for double super Yangians has not yet been investigated. In this paper, we study the type $\boldsymbol{A}$ double super Yangians associated with arbitrary sequences $\mathfrak{s}$.

The main results of this paper are threefold. First, we define the double Yangian of the Lie superalgebras $\mathfrak{gl}_{m|n}^{\mathfrak{s}}$ and $\mathfrak{sl}_{m|n}^{\mathfrak{s}}$ associated with an arbitrary sequence $\mathfrak{s}$, i.e., double super Yangians $\mathrm{DY}_{h}(\mathfrak{gl}_{m|n}^{\mathfrak{s}})$ and $\mathrm{DY}_{h}(\mathfrak{sl}_{m|n}^{\mathfrak{s}})$, via the R-matrix presentation. Next, we present their Drinfeld presentations using the Gauss decomposition of the generator matrix and establish the isomorphism between the two constructions following the approach of \cite{YJ24}. Finally, we construct level-1 representations of these double super Yangians in terms of bosons. To this end,  we construct a Poincar\'e-Birkhoff-Witt type basis for $\mathrm{DY}_{h}(\mathfrak{gl}_{m|n}^{\mathfrak{s}})$ and prove that it specializes to the universal enveloping algebra of the affine Lie superalgebra $\widehat{\mathfrak{gl}}^{\mathfrak{s}}_{m|n}$. These two key tools are necessary for proving the equivalence of the two presentations. We also generalize the quantum Berezinian for $\mathrm{DY}_{h}(\mathfrak{gl}_{m|n}^{\mathfrak{s}})$ from standard parity to arbitrary cases. Following the approach of \cite{XLZ}, this quantum Berezinian can be used to define the R-matrix presentation of $\mathrm{DY}_{h}(\mathfrak{sl}_{m|n}^{\mathfrak{s}})$.

This paper is organized as follows. Section 2 collects the notation and preliminaries used throughout the paper. In Section 3, the double super Yangian $\mathrm{DY}_{h}(\mathfrak{gl}_{m|n}^{\mathfrak{s}})$ is defined via the R-matrix presentation, and its Poincar{\'e}-Birkhoff-Witt theorem is established. Section 4 is devoted to deriving the relations for the Gaussian generators. The isomorphism between the Drinfeld and R-matrix presentations of $\mathrm{DY}_{h}(\mathfrak{gl}_{m|n}^{\mathfrak{s}})$ is established in Section 5. Section 6 provides the quantum Berezinian for $\mathrm{DY}_{h}(\mathfrak{gl}_{m|n}^{\mathfrak{s}})$. In Section 7, we introduce both the R-matrix and the Drinfeld presentations of $\mathrm{DY}_{h}(\mathfrak{sl}_{m|n}^{\mathfrak{s}})$ and prove the equivalence between them. Finally, Section 8 presents a bosonic construction of level 1 modules for $\mathrm{DY}_{h}(\mathfrak{gl}_{m|n}^{\mathfrak{s}})$ and $\mathrm{DY}_{h}(\mathfrak{sl}_{m|n}^{\mathfrak{s}})$.

\medskip
\section{Preliminaries}
In this section, we recall the notations and basic definitions that will be used later. We adopt the following conventions: let $\mathbb{C}$ be the complex numbers, $\mathbb{C}^*$ the nonzero complex numbers, $\mathbb{Z}$ the integers, $\mathbb{N}$ the non-negative integers, $\mathbb{Z}^*$ the nonzero integers and $\mathbb{Z}_+$ the positive integers. We write $\mathbb{Z}_2=\mathbb{Z}/2\mathbb{Z}=\{\bar{0},\bar{1}\}$.  Let $\delta_{ij}$ be the Kronecker symbol, which equals 1 if $i=j$ and 0 otherwise.

All superspaces, (associative) superalgebras, and Lie superalgebras considered in this section are defined over $\mathbb{C}$.
A \textit{superspace} is a $\mathbb{Z}_2$-graded vector space $V=V_{\bar{0}} \oplus V_{\bar{1}}$. Elements of $V_{\bar{0}}$ are called even, and elements of $V_{\bar{1}}$ are called odd. An element $x\in V$ is  \textit{homogeneous} if it is even or odd. For a homogeneous element $x\in V$, its parity $|x|$ is defined to be 0 for even elements and 1 for odd elements, i.e., $|x|=0$ when $x\in V_{\bar{0}}$ and $|x|=1$ when $x\in V_{\bar{1}}$.
An (associative) \textit{superalgebra} is a $\mathbb{Z}_2$-graded algebra $U=U_{\bar{0}}\oplus U_{\bar{1}}$ such that $xy\in U_{\bar{i}+\bar{j}}$ for homogeneous elements $x\in U_{\bar{i}}$, $y\in U_{\bar{j}}$ with $i,j\in\mathbb{Z}_2$. For two such superalgebras $U$ and $W$, the tensor product $U\otimes W$ admits a natural superalgebra structure, satisfying the following multiplication rule
\begin{gather*}
    (x_1 \otimes y_1)(x_2 \otimes y_2)=(-1)^{|y_1||x_2|}(x_1x_2 \otimes y_1y_2),
\end{gather*}
where $x_1,x_2\in U$ and  $y_1,y_2\in W$ are homogeneous elements. In particular, the even and odd parts of $(U\otimes W)$ are given by
\begin{align*}
    (U\otimes W)_{\bar{0}}&=(U_{\bar{0}}\otimes W_{\bar{0}})\oplus (U_{\bar{1}}\otimes W_{\bar{1}}), \\
    (U\otimes W)_{\bar{1}}&=(U_{\bar{0}}\otimes W_{\bar{1}})\oplus (U_{\bar{1}}\otimes W_{\bar{0}}).
\end{align*}

A \textit{Lie superalgebra} is a superspace $G=G_{\bar{0}} \oplus G_{\bar{1}}$ equipped with a bilinear map $[\,\cdot,\, \cdot\,]: G\times G\to G$, called the \textit{Lie superbracket}, satisfying the following axioms for homogeneous elements $x,y,z\in G$:
\begin{align*}
    &[x,\,y]=-(-1)^{|x||y|}[y,\,x],  &&(\text{skew-supersymmetry}) \\
    &(-1)^{|x||y|}[[x,\,y],\,z]+(-1)^{|y||z|}[[y,\,z],\,x]+(-1)^{|z||x|}[[z,\,x],\,y]=0  &&(\text{super-Jacobi identity})
\end{align*}
Subsequently, we review and restate the definition and properties of Lie superalgebras of type $\boldsymbol{A}$. This follows \cite{K77,LZ25,M22,Mu12} and the references therein.

Given $m,n\in\mathbb{N}$ with $N=m+n\geqslant 2$. Let $\mathcal{S}_{m|n}$ be the set of all $0^{m}1^{n}$-sequences $\mathfrak{s}=(s_1s_2\cdots s_N)$ consisting of $m$ zeros and $n$ ones. Each $\mathfrak{s}\in\mathcal{S}_{m|n}$ is called a \textit{parity sequence}. The \textit{standard parity} sequence, denoted by $\mathfrak{s}^{st}$, is the one with $s_i=0$ for $i=1,\cdots,m$, and $s_i=1$ for $i=m+1,\cdots,N$.
Let $I_{\mathfrak{s}}=\{1,\ldots,N\}$ be the index set. For a fixed $\mathfrak{s}\in\mathcal{S}_{m|n}$, we define two functions on $I_{\mathfrak{s}}$ by
$$
 |i|=\begin{cases}
        0, &\text{if}~~s_i=0; \\
        1, &{\rm otherwise}.
    \end{cases}
    \qquad
    d_i=\begin{cases}
        1, &\text{if}~~s_i=0; \\
        -1, &{\rm otherwise}.
    \end{cases}
$$
Throughout this paper, $I_{\mathfrak{s}}$ will always denote the index set equipped with these two functions. Whenever we write $x_{ij}$ for $i,j\in I_{\mathfrak{s}}$, it is understood that $x_{ij}$ is homogeneous of parity $|i|+|j|$ (with respect to $\mathfrak{s}$); the definition extends to arbitrary elements by linearity.

\begin{rema}\label{conven-R}
We should emphasize here that, unless otherwise specified, the same notation is used for the generators of the (Lie) superalgebras for different $\mathfrak{s}\in\mathcal{S}_{m|n}$, with the understanding that the relevant superalgebra will always be clear from the context.
\end{rema}

Fix $\mathfrak{s}\in\mathcal{S}_{m|n}$. Let $e_1,\ldots,e_N$ denote the standard basis of the $m|n$-dimensional superspace $\mathcal{V}_{\mathfrak{s}}=\mathbb{C}^{m|n}$ with parity $|e_{i}|=|i|$. The elementary matrices $E_{ij}$ for $i,j\in I_{\mathfrak{s}}$ constitute a basis for the endomorphism ring $\operatorname{End} \mathcal{V}_{\mathfrak{s}}$ as a linear space, satisfying
$$
E_{ij}(e_{k})=\delta_{jk}e_{i},\quad~\forall\, i,j,k\in I_{\mathfrak{s}}.
$$

Let $U$ be a superalgebra. Given a $\mathbb{Z}_2$-graded matrix
$X=(x_{ij})_{i,j\in I_{\mathfrak{s}}}=\sum_{i,j\in I_{\mathfrak{s}}} x_{ij} \otimes E_{ij}  \in U\otimes \operatorname{End}\mathcal{V}_{\mathfrak{s}}$ and an integer $t\geqslant 2$.
For each $1\leqslant a\leqslant t$, we denote by $X_a$ its $a$-th copy of the superspace $\mathrm{End}\mathcal{V}_{\mathfrak{s}}$ by
$$
X_a=\sum\limits_{i,j\in I_{\mathfrak{s}}}x_{ij}\otimes 1^{\otimes (a-1)} \otimes E_{ij}\otimes 1^{\otimes (t-a)} \in U\otimes \mathrm{End}\mathcal{V}_{\mathfrak{s}}^{\otimes t}.
$$
For $R=\sum_i r_{i}\otimes r^{i}\in \operatorname{End}\mathcal{V}_{\mathfrak{s}}^{\otimes 2}$, define the element $R_{ab} \in \operatorname{End}\mathcal{V}_{\mathfrak{s}}^{\otimes t}$ for $1\leqslant a<b\leqslant t$ by
\begin{gather*}
    R_{ab}=\sum_{i} 1^{\otimes (a-1)}\otimes r_i\otimes 1^{\otimes (b-a-1)}\otimes r^i\otimes 1^{\otimes (t-b)}.
\end{gather*}
For example, when $t=3$, we have
$$
R_{12}=R\otimes 1, \quad R_{13}=\sum_{i} r_{i}\otimes 1 \otimes r^{i}, \quad R_{23}=1\otimes R.
$$
Let $P$ be the $\mathbb{Z}_2$-graded permutation operator on the tensor product $\mathcal{V}_{\mathfrak{s}}^{\otimes 2}$ given by
$$
P=\sum_{i,j\in I_{\mathfrak{s}}}d_{j}E_{ij}\otimes E_{ji}.
$$
We then define $R_{ba}=P_{ab} R_{ab}P_{ab}$.

The endomorphism ring $\operatorname{End}\mathcal{V}_{\mathfrak{s}}$ admits a Lie superalgebra structure with the superbracket
\begin{align*}
\left[E_{ij},E_{kl}\right]=\delta_{jk}E_{il}-\epsilon_{ij;kl}\delta_{il}E_{kj},
\end{align*}
where $\epsilon_{ij;kl}=(-1)^{(|i|+|j|)(|k|+|l|)}$.
 Equipped with this bracket, $\operatorname{End}\mathcal{V}_{\mathfrak{s}}$ is called the \textit{general linear Lie superalgebra}, denoted by $\mathfrak{gl}_{m|n}^{\mathfrak{s}}$.

The \textit{supertrace} map $\mathbf{str}$ is a $\mathbb{C}$-valued linear functional on $\mathfrak{gl}_{m|n}^{\mathfrak{s}}$, determined by $E_{ij}\mapsto d_{i}\delta_{ij}$ for $i,j\in I_{\mathfrak{s}}$. Then the \textit{special linear Lie superalgebra} $\mathfrak{sl}_{m|n}^{\mathfrak{s}}$ is defined as
$$
\mathfrak{sl}_{m|n}^{\mathfrak{s}}:=\{x\in\mathfrak{gl}_{m|n}^{\mathfrak{s}}|\mathbf{str}(x)=0\}.
$$
Let $\mathcal{I}=\sum_{i\in I_{\mathfrak{s}}}E_{ii}$ be the identity matrix.
We note that the Lie superalgebra $\mathfrak{sl}_{m|n}^{\mathfrak{s}}$ is simple unless $m=n$, in which case it has a one-dimensional ideal $\langle \mathcal{I}\rangle$ spanned by scalar matrices $\lambda \mathcal{I}$, $\lambda \in\mathbb{C}$.

Following \cite{K77}, the family $\boldsymbol{A}^{\mathfrak{s}}(m,n)$ of classical simple Lie superalgebras is given by:
\begin{align*}
\boldsymbol{A}^{\mathfrak{s}}(m-1,n-1)&=\mathfrak{sl}_{m|n}^{\mathfrak{s}}\qquad~~\text{for}~~m\neq n,\quad m,n\geqslant 1;\\
\boldsymbol{A}^{\mathfrak{s}}(n-1,n-1)&=\mathfrak{sl}_{n|n}^{\mathfrak{s}}/\langle \mathcal{I}\rangle := \mathfrak{psl}^{\mathfrak{s}}_{n|n},\quad \text{for} ~~n>1.
\end{align*}
A Lie superalgebra $\mathfrak{g}$ is said to be of type $\boldsymbol{A}$ if $\mathfrak{g}=\mathfrak{gl}_{m|n}^{\mathfrak{s}}$, $\mathfrak{sl}_{m|n}^{\mathfrak{s}}$ or $\mathfrak{psl}^{\mathfrak{s}}_{n|n}$.

Let $\mathfrak{h}_{\mathfrak{s}}:=\text{Span}_{\mathbb{C}}\{h_{i}:=d_{i}E_{ii}\ |\ i\in I_{\mathfrak{s}}\}$ be the \textit{Cartan subalgebra} of $\mathfrak{gl}_{m|n}^{\mathfrak{s}}$. The linear functions $\epsilon_{i}\in\mathfrak{h}_{\mathfrak{s}}^{\ast}$ are defined by $\epsilon_{i}(h_{j})=d_{i}\delta_{ij}$, and their parities are given by $|\epsilon_{i}|=|i|$. Set $I_{\mathfrak{s}}'=I_{\mathfrak{s}}\setminus\{N\}$. For each $i\in I_{\mathfrak{s}}'$, we denote
$\alpha_{i}=\epsilon_{i}-\epsilon_{i+1}$  to be the \textit{simple root} and $h_{\alpha_{i}}=h_{i}-h_{i+1}$ the \textit{simple coroot} of type $\boldsymbol{A}$. We call $\mathcal{Q}_{\mathfrak{s}}=\sum_{i\in I_{
\mathfrak{s}}'}\mathbb{Z}\alpha_{i}$  the \textit{root lattice} and $\mathcal{Q}_{\mathfrak{s}}^{\vee}=\sum_{i\in I_{\mathfrak{s}}'}\mathbb{Z}h_{\alpha_{i}}$ the \textit{coroot lattice}. The Cartan matrix of the Lie superalgebra in type $\boldsymbol{A}$ is defined as the $(N-1)\times (N-1)$-matrix $A=\left(a_{i j}\right)_{i,j\in I_{\mathfrak{s}}'}$ with entries
\begin{gather*}
    a_{i j}=\left(d_{i}+d_{i+1}\right) \delta_{i, j}-d_{i}\delta_{i, j+1}-d_{i+1}\delta_{i+1, j}.
\end{gather*}

Let $\mathfrak{g}$ be the Lie superalgebra of type $\boldsymbol{A}$. The affine Lie superalgebra $\widehat{\mathfrak{g}}$ is defined as a one-dimensional central extension of the loop superalgebra associated with $\mathfrak{g}$. That is, $\widehat{\mathfrak{g}}=\mathfrak{g}
\otimes \mathbb{C}[t,t^{-1}]\oplus \mathbb{C}K$, where the element $K$ is central.
The superbracket on $\widehat{\mathfrak{g}}$ is defined by
\begin{align}\label{affine}
\left[E_{ij}^{(r)},E_{kl}^{(s)}\right]=\delta_{kj}E_{il}^{(r+s)}-\delta_{il}\epsilon_{ij;kl}
E_{kj}^{(r+s)}+rd_{i}\delta_{kj}\delta_{il}\delta_{r+s,0}K,
\end{align}
where $E_{ij}^{(r)}:=E_{ij}\otimes t^r\in\widehat{\mathfrak{g}}$ for $i,j\in I_{\mathfrak{s}}$ and $r\in \mathbb{Z}$.

\medskip
\section{Double super Yangian of the general linear Lie superalgebra}

In this section, we define the double Yangian associated with the general linear Lie superalgebra $\mathfrak{gl}^{\mathfrak{s}}_{m|n}$ in terms of R-matrix generators, for any parity sequence $\mathfrak{s}\in\mathcal{S}_{m|n}$. We call it \textit{double super Yangian} (of $\mathfrak{gl}^{\mathfrak{s}}_{m|n}$) and denote it by $\mathrm{DY}_{h}(\mathfrak{gl}_{m|n}^{\mathfrak{s}})$. Moreover, we prove that, for any given $(m,n)$, the resulting definition is independent of the chosen sequence $\mathfrak{s}$, and we then establish the PBW basis of $\mathrm{DY}_{h}(\mathfrak{gl}_{m|n}^{\mathfrak{s}})$ for each $\mathfrak{s}$. Notice that the definition of the standard version $\mathrm{DY}_{-h}(\mathfrak{gl}_{m|n}^{\mathfrak{s}^{st}})$ and its PBW basis were given in \cite{BK25-1, BK25-2}.

\subsection{R-matrix presentation}

Let $h$ be a formal parameter, and $\mathcal{A}:=\mathbb{C}[[h]]$ the ring of formal power series endowed with the $h$-adic topology. For a given $\mathfrak{s}\in\mathcal{S}_{m|n}$, the rational R-matrix $R(u)$ is defined by
\begin{equation}\label{R-matrix}
R(u)=1+u^{-1}hP\in \operatorname{End} \mathcal{V}_{\mathfrak{s}}^{\otimes 2}\otimes[u^{-1}],
\end{equation}
where $R(u)$ possesses an expansion in negative powers of $u$ or positive powers of $h$. It satisfies the Yang–Baxter equation and the unitary condition with values in $\operatorname{End} \mathcal{V}_{\mathfrak{s}}^{\otimes3}$ :
\begin{align}\notag
R_{12}(u-v)R_{13}(u)R_{23}(v)&=R_{23}(v)R_{13}(u)R_{12}(u-v),\\ \label{R12R21}
R_{12}(u)R_{21}(-u)&=1-h^{2}u^{-2}.
\end{align}

\begin{defi}\label{DY}
For a given $\mathfrak{s}\in\mathcal{S}_{m|n}$, the double super Yangian $\mathrm{DY}_{h}(\mathfrak{gl}_{m|n}^{\mathfrak{s}})$ is defined as the $\mathbb{Z}_{2}$-graded unital associative superalgebra over $\mathcal{A}$ generated by the elements $\{t_{ij}^{(\pm r)}|i,j\in I_{\mathfrak{s}}, r\in \mathbb{Z}_{+}\}$ and an even central element $C$. The defining relations are expressed in terms of the generating series
\begin{align*}
t_{ij}^{+}(u)=\delta_{ij}-h\sum_{r\geqslant 1}t_{ij}^{(r)}u^{-r}, \qquad
t_{ij}^{-}(u)=\delta_{ij}+h\sum_{r\geqslant 1}t_{ij}^{(-r)}u^{r-1}.
\end{align*}
Set $T^{\pm}(u)=\left(t_{ij}^{\pm}(u)\right)_{i,j\in I_{\mathfrak{s}}}$. The defining relations are given as follows$\mathrm{:}$
\begin{align}\label{RTT1}
R_{12}(u-v)T_{1}^{\pm}(u)T_{2}^{\pm}(v)&=T_{2}^{\pm}(v)T_{1}^{\pm}(u)R_{12}(u-v),\\ \label{RTT2}
R_{12}(u-v-\frac{1}{2}hC)T_{1}^{+}(u)T_{2}^{-}(v)&=T_{2}^{-}(v)T_{1}^{+}(u)R_{12}(u-v+\frac{1}{2}hC),
\end{align}
where the R-matrix $R(u)$ is viewed as power series in $u^{-1}$.
\end{defi}

In the special case $n=0$, the double Yangian $\mathrm{DY}_{h}(\mathfrak{gl}_{m|n}^{\mathfrak{s}})$ descends to the double Yangian $\mathrm{DY}_{h}(\mathfrak{gl}_{m})$ studied by Yang and Jing \cite{YJ24}.
In analogy to \cite{BK25-2}, the superalgebra $\mathrm{DY}_{h}(\mathfrak{gl}_{m|n}^{\mathfrak{s}})$ admits a natural Hopf superalgebra structure defined by
\begin{align*}
&\Delta(t_{ij}^{\pm}(u))=\sum_{k=1}^{N}t_{ik}^{\pm}(u\mp\frac{1}{4}hC_2)\otimes t_{kj}^{\pm}(u\pm\frac{1}{4}hC_1),\quad \Delta(C)=C_1+C_2, \\
&\varepsilon(T^{\pm}(u))=1,~~\varepsilon(C)=0,~~S(T^{\pm}(u))=T^{\pm}(u)^{-1},~~S(C)=-C, 
\end{align*}
where $C_{1}=C\otimes 1$ and $C_{2}=1\otimes C$.

We define the level-0 double Yangian of $\mathfrak{gl}_{m|n}^{\mathfrak{s}}$ by $\mathrm{DY}_{h}^{0}(\mathfrak{gl}_{m|n}^{\mathfrak{s}}):=\mathrm{DY}_{h}(\mathfrak{gl}_{m|n}^{\mathfrak{s}})/\langle C\rangle$, that is, by setting the central element $C$ to zero. The algebra $\mathrm{DY}_{h}(\mathfrak{gl}_{m|n}^{\mathfrak{s}})$ is then a central extension of $\mathrm{DY}_{h}^{0}(\mathfrak{gl}_{m|n}^{\mathfrak{s}})$ obtained by adjoining the central generator $C$, and there is an isomorphism
\begin{align*}
 \mathrm{DY}_{h}(\mathfrak{gl}_{m|n}^{\mathfrak{s}})\cong \mathrm{DY}_{h}^{0}(\mathfrak{gl}_{m|n}^{\mathfrak{s}})\otimes_{\mathcal{A}}\mathcal{A}[C].   
\end{align*}

\begin{rema}
The subalgebra over $\mathbb{C}[h]$ generated by all elements $t_{ij}^{(r)}$ $($resp. $t_{ij}^{(-r)})$ with $i,j\in I_{\mathfrak{s}}$ and $r\in\mathbb{Z}_{+}$ is precisely the super Yangian $\mathrm{Y}^{+}_{h}(\mathfrak{gl}_{m|n}^{\mathfrak{s}})$ $($resp. dual super Yangian $\mathrm{Y}_{h}^{-}(\mathfrak{gl}_{m|n}^{\mathfrak{s}}))$. Throughout this paper, we assume that all these superalgebras are complete with respect to the $h$-adic topology.
\end{rema}

The defining relations \eqref{RTT1} and \eqref{RTT2} give rise to the following commutation relations:
\begin{align}\label{Y1}
\left[t_{ij}^{+}(u),t_{kl}^{+}(v)\right]&=\theta_{i,j,k}\frac{h}{v-u}\left(t_{kj}^{+}(u)t_{il}^{+}(v)-t_{kj}^{+}(v)t_{il}^{+}(u)\right),\\ \label{Y2}
\left[t_{ij}^{-}(u),t_{kl}^{-}(v)\right]&=\theta_{i,j,k}\frac{h}{v-u}\left(t_{kj}^{-}(u)t_{il}^{-}(v)-t_{kj}^{-}(v)t_{il}^{-}(u)\right),\\ \label{Y3}
\left[t_{ij}^{+}(u),t_{kl}^{-}(v)\right]&=\theta_{i,j,k}\left(
\frac{h}{v_{-}-u_{+}}t_{kj}^{+}(u)t_{il}^{-}(v)-\frac{h}{v_{+}-u_{-}}t_{kj}^{-}(v)t_{il}^{+}(u)\right),
\end{align}
where $\theta_{i,j,k}=(-1)^{|i||j|+|i||k|+
|j||k|}$ and  $u_{\pm}=u\pm\frac{1}{4}hC$.

Using the definition of the R-matrix together with the unitary relation \eqref{R12R21}, the defining relations \eqref{RTT1} and \eqref{RTT2} imply the following lemma.

\begin{lemm}
    The following relations hold in $\mathrm{DY}_{h}(\mathfrak{gl}_{m|n}^{\mathfrak{s}})\otimes
(\operatorname{End}\mathcal{V}_{\mathfrak{s}}^{\otimes 2})[[u^{\pm 1},v^{\pm 1}]]:$
    \begin{align}\label{RTT4}
&T_{1}^{\pm}(v)^{-1}R_{21}(u-v) T_{2}^{\pm}(u)=T_{2}^{\pm}(u)R_{21}(u-v)T_{1}^{\pm}(v)^{-1},\\ \label{RTT5}
&T_{1}^{+}(v)^{-1}R_{21}(u_{-}-v_{+}) T_{2}^{-}(u)=T_{2}^{-}(u)R_{21}(u_{+}-v_{-})T_{1}^{+}(v)^{-1},\\ \label{RTT6}
&T_{2}^{\pm}(u)^{-1}T_{1}^{\pm}(v)^{-1}R_{21}(u-v)=
R_{21}(u-v)T_{1}^{\pm}(v)^{-1}T_{2}^{\pm}(u)^{-1},\\ \label{RTT7}
&T_{2}^{-}(u)^{-1}T_{1}^{+}(v)^{-1}R_{21}(u_{-}-v_{+})
=R_{21}(u_{+}-v_{-})T_{1}^{+}(v)^{-1}T_{2}^{-}(u)^{-1},\\ \label{RTT8}
&\frac{(u_{-}-v_{+})^{2}}{(u_{-}-v_{+})^{2}-h^{2}}T_{1}^{-}(v)T_{2}^{+}(u) R_{21}(u_{-}-v_{+})
=\frac{(u_{+}-v_{-})^{2}}{(u_{+}-v_{-})^{2}-h^{2}}R_{21}(u_{+}-v_{-})
T_{2}^{+}(u)T_{1}^{-}(v),\\  \label{RTT9}
&\frac{(u_{-}-v_{+})^{2}}{(u_{-}-v_{+})^{2}-h^{2}}T_{2}^{+}(u)R_{21}(u_{-}-v_{+})T_{1}^{-}(v)^{-1}
=\frac{(u_{+}-v_{-})^{2}}{(u_{+}-v_{-})^{2}-h^{2}}T_{1}^{-}(v)^{-1}R_{21}(u_{+}-v_{-})T_{2}^{+}(u),
\end{align}
\end{lemm}

Similarly to \cite[Lemma 2.24]{T20}, the following result shows that, for fixed $(m,n)$, the above definition of $\mathrm{DY}_{h}(\mathfrak{gl}_{m|n}^{\mathfrak{s}})$ is independent of the choice of $\mathfrak{s}\in\mathcal{S}_{m|n}$.

\begin{lemm}\label{DY-iso}
Up to isomorphism, the double super Yangian $\mathrm{DY}_{h}(\mathfrak{gl}_{m|n}^{\mathfrak{s}})$ is independent of the parity sequence $\mathfrak{s}$ and depends only on $(m,n)$.
\end{lemm}
\begin{proof}
Let $\mathfrak{s}_1$ and $\mathfrak{s}_2$ be two arbitrary elements of $\mathcal{S}_{m|n}$, and let $\mathcal{V}_{\mathfrak{s}_1}$ and $\mathcal{V}_{\mathfrak{s}_2}$ be the corresponding superspaces, equipped with $\mathbb{C}$-bases $\{v_{i}\}_{i\in I_{\mathfrak{s}_1}}$ and $\{v_{i}'\}_{i\in I_{\mathfrak{s}_2}}$, respectively. Observe that the dimensions of even and odd subspaces satisfy
$$
\dim((\mathcal{V}_{\mathfrak{s}_1})_{\bar{0}})=\dim((\mathcal{V}_{\mathfrak{s}_2})_{\bar{0}})=m,\quad \dim((\mathcal{V}_{\mathfrak{s}_1})_{\bar{1}})=\dim((\mathcal{V}_{\mathfrak{s}_2})_{\bar{1}})=n.
$$
Then we may choose a permutation $\sigma$ of order $N$ such that for all $i\in I_{\mathfrak{s}_1}$, we have $\sigma(i)\in I_{\mathfrak{s}_2}$ and $|v_i|=|v_{\sigma(i)}'|$.
Moreover, the assignment $t_{ij}^{(\pm r)}\mapsto t_{\sigma(i),\sigma(j)}^{(\pm r)}$ preserves the defining relations \eqref{Y1}-\eqref{Y3} and is invertible, thereby yielding an isomorphism $\mathrm{DY}_{h}(\mathfrak{gl}_{m|n}^{\mathfrak{s}_1})\rightarrow \mathrm{DY}_{h}(\mathfrak{gl}_{m|n}^{\mathfrak{s}_2})$.
\end{proof}

\subsection{Poincar{\'e}-Birkhoff-Witt type basis}
For simplicity, we consider the double super Yangians at level-0 first. We set $t_{ij}^{(0)}: =\delta_{ij}$.
In $\mathrm{DY}_{h}^{0}(\mathfrak{gl}^{\mathfrak{s}}_{m|n})$, the relation \eqref{RTT2} implies the following commutation relation between the generators $t_{ij}^{(r)}$ and $t_{kl}^{(-s)}$:
\begin{equation}\label{t1}
    \begin{split}
        \left[t_{ij}^{(r)},t_{kl}^{(-s)}\right]=&\theta_{i,j,k}\left(\delta_{s < r}\left(\delta_{il}t_{kj}^{(r-s)}-\delta_{kj}t_{il}^{(r-s)}\right)-\delta_{s \geq r}\left(\delta_{il}t_{kj}^{(r-s-1)}-\delta_{kj}t_{il}^{(r-s-1)}\right)\right)\\ 
+&\theta_{i,j,k}h\sum_{p=1}^{r-1}\left(t_{kj}^{(p)}t_{il}^{(r-s-p-1)}-
t_{kj}^{(r-s-p-1)}t_{il}^{(p)}\right).
    \end{split}
\end{equation}

\begin{lemm}
In $\mathrm{DY}_{h}^{0}(\mathfrak{gl}_{m|n}^{\mathfrak{s}})\otimes
(\operatorname{End} \mathcal{V}_{\mathfrak{s}}^{\otimes 2})[[u^{\pm 1},v^{\pm 1}]]$, the following equation holds:
\begin{align}\label{RTT3}
R_{12}(u-v)T_{1}^{+}(u)T_{2}^{-}(v)=T_{2}^{-}(v)T_{1}^{+}(u)R_{12}(u-v).
\end{align}
\end{lemm}
\begin{proof}
From \eqref{RTT2}, multiplying on both sides by $R_{21}(v-u)$ and using \eqref{R12R21} gives
\begin{align*}
T_{1}^{-}(u)T_{2}^{+}(v)R_{21}(v-u)=R_{21}(v-u)T_{2}^{+}(v)T_{1}^{-}(u).
\end{align*}
After applying the operator $P_{12}(\cdot)P_{12}$ to this equation and interchanging $u$ and $v$, \eqref{RTT3} follows from $R_{21}=P_{12}R_{12}P_{12}$.
\end{proof}

From \eqref{RTT3}, we also obtain the following commutation relation:
\begin{equation}\label{t2}
    \begin{split}
        \left[t_{kl}^{(-s)},t_{ij}^{(r)}\right]=&\theta_{k,l,i}\left(\delta_{s \geq r}\left(\delta_{il}t_{kj}^{(r-s-1)}-\delta_{kj}t_{il}^{(r-s-1)}\right)+
        \delta_{s < r}\left(\delta_{il}t_{kj}^{(r-s)}-\delta_{kj}t_{il}^{(r-s)}\right)\right)\\ 
        +&\theta_{k,l,i}h\sum_{p=1}^{r-1}\left(t_{il}^{(r-s-p-1)}t_{kj}^{(p)}-t_{il}^{(p)}t_{kj}^{(r-s-p-1)}\right).
    \end{split}
\end{equation}

Let $\prec$ be the lexicographical order on the countable set $I_{\mathfrak{s}}\times I_{\mathfrak{s}} \times \mathbb{Z}_{+}$. For any $i_1,i_2,j_1,j_2 \in I_{\mathfrak{s}},r_1,r_2 \in \mathbb{Z}_{+}$, we write $(i_1,j_1,r_1)\prec(i_2,j_2,r_2)$ if and only if one of the following conditions holds:
$$
(1)~i_1<i_2,\quad(2)~i_1=i_2,~ j_1<j_2,\quad (3)~i_1=i_2,~ j_1=j_2,~r_1<r_2.
$$
Using this order, we define an ordering on the generators of $\mathrm{DY}_{h}^{0}(\mathfrak{gl}^{\mathfrak{s}}_{m|n})$ as follows:
\begin{align*}
&(i)~\gamma_{i_1,j_1}^{(r_1)}\prec \gamma_{i_2,j_2}^{(r_2)} \quad {\rm if~and~only~if}\quad (j_1-i_1,i_1,r_1)\prec (j_2-i_2,i_2,r_2),\\
&(ii)~t_{ij}^{(-r)}\prec t_{ij}^{(r)}\quad {\rm for~any~triples}\quad (i,j,r),
\end{align*}
where $\gamma_{i,j}^{(r)}\in\{t_{i,j}^{(r)},t_{i,j}^{(-r)}\}$. The following theorem establishes a PBW basis for $\mathrm{DY}_{h}^{0}(\mathfrak{gl}^{\mathfrak{s}}_{m|n})$ with respect to this ordering.

\begin{lemm}
Let $\mathfrak{B}$ denote the set of all monomials in the generators $t_{ij}^{(\pm r)}$ ($i,j\in I_{\mathfrak{s}}$, $r\in\mathbb{Z}_+$), with the generators ordered by $\prec$ and with odd generators appearing with exponents at most 1.
Then $\mathfrak{B}$ forms a basis of $\mathrm{DY}_{h}^{0}(\mathfrak{gl}^{\mathfrak{s}}_{m|n})$ over $\mathcal{A}$.
\end{lemm}
\begin{proof}
We first show that $\mathfrak{B}$ spans $\mathrm{DY}_{h}^{0}(\mathfrak{gl}^{\mathfrak{s}}_{m|n})$. We only need to show 
$$t_{ij}^{(r)}t_{kl}^{(s)},~~  t_{ij}^{(r)}t_{kl}^{(-s)},~~ t_{ij}^{(-r)}t_{kl}^{(s)},~~ t_{ij}^{(-r)}t_{kl}^{(-s)}\in\operatorname{Span}_{\mathcal{A}}\mathfrak{B}.$$
From \eqref{t1} and \eqref{t2}, we get
\begin{align}\label{t3}
&t_{ij}^{(r)}t_{kl}^{(-s)}\in\operatorname{Span}_{\mathcal{A}}\{t_{kl}^{(-s)}t_{ij}^{(r)}, t_{kj}^{(p)}t_{il}^{(r-s-p-1)}, t_{kj}^{(r-s-p-1)}t_{il}^{(p)}, t_{kj}^{(r-s)},t_{il}^{(r-s)}\mid 1\leq p\leq r-1\},\\ \label{t4}
&t_{ij}^{(r)}t_{kl}^{(-s)}\in\operatorname{Span}_{\mathcal{A}}\{t_{kl}^{(-s)}t_{ij}^{(r)}, t_{il}^{(p)}t_{kj}^{(r-s-p-1)}, t_{il}^{(r-s-p-1)}t_{kj}^{(p)}, t_{kj}^{(r-s)},t_{il}^{(r-s)}\mid 1\leq p\leq r-1\}.
\end{align}
Suppose $j-i>l-k$, which implies that $t_{kl}^{(-s)}$ precedes $t_{ij}^{(r)}$. Regardless of whether $j-k<l-i$ or $j-k>l-i$, it follows from \eqref{t3} and \eqref{t4} that
\begin{align}\label{t5}
t_{ij}^{(r)}t_{kl}^{(-s)}\in\mathfrak{B}.
\end{align}
In the remaining subcase $j-k=l-i$, we have
$$
j-k=l-i\Rightarrow j+i=k+l\Rightarrow j-i=2(k-i)+l-k,
$$
hence $i<k$. According to \eqref{t4}, we still deduce \eqref{t5}.

Now consider the case $j-i=l-k$. If $i>k$, then $j-k=l-i+2(i-k)>l-i$, so \eqref{t5} follows from \eqref{t3}. If $i=k$ and $r>s$, we have $s-r+p+1=p+1-(r-s)\leq p$. By \eqref{t2}, we have $t_{kl}^{(-s)}t_{ij}^{(r)}$, $t_{il}^{(r-s-p-1)}t_{kj}^{(p)}
\in\mathfrak{B}$ and
\begin{align*}
t_{il}^{(p)}t_{kj}^{(r-s-p-1)}\in\operatorname{Span}_{\mathbb{C}}\{t_{kj}^{(r-s-p-1)}t_{il}^{(p)}, t_{ij}^{(a)}t_{kl}^{(r-s-a-2)}, t_{ij}^{(r-s-a-2)}t_{kl}^{(a)}\mid 1\leq a\leq r-2\}.
\end{align*}
Among these, $t_{kj}^{(r-s-p-1)}t_{il}^{(p)}\in\mathfrak{B}$, and exactly one of $t_{ij}^{(a)}t_{kl}^{(r-s-a-2)}\in\mathfrak{B}$ and $t_{ij}^{(r-s-a-2)}t_{kl}^{(a)}\in\mathfrak{B}$. For the term not in $\mathfrak{B}$, we apply the argument recursively. Since each application decreases the relevant index, the process terminates after finitely many steps. So, we can deduce that $t_{il}^{(p)}t_{kl}^{(r-s-p-1)}\in\mathfrak{B}$. Hence, \eqref{t5} holds.

By the same argument, we obtain $t_{ij}^{(r)}t_{kl}^{(s)},~ t_{ij}^{(-r)}t_{kl}^{(s)}, ~t_{ij}^{(-r)}t_{kl}^{(-s)}\in\operatorname{Span}_{\mathcal{A}}\mathfrak{B}.
$

Finally, since the proof of linear independence in \cite[Theorem 2.6]{BK25-2} does not depend on the generator ordering or the parity sequence, the same argument establishes that $\mathfrak{B}$ is linearly independent.
\end{proof}

Extend the total order $(i)$ and $(ii)$ to $\mathrm{DY}_{h}(\mathfrak{gl}^{\mathfrak{s}}_{m|n})$ by inserting the central element $C$ in an arbitrary manner. Arguing as in the proof of \cite[Theorem 2.9]{BK25-2}, we obtain a PBW basis for $\mathrm{DY}_{h}(\mathfrak{gl}^{\mathfrak{s}}_{m|n})$ with respect to this extended order.

\begin{theo}\label{PBW basis}
Let $\mathfrak{B}^C$ denote the set of all monomials
in the generators $C$ and $t_{ij}^{(\pm r)}$ ($i,j\in I_{\mathfrak{s}}$, $r\in\mathbb{Z}_+$) ordered by $\prec$, in which the odd generators appear with exponents at most 1. Then $\mathfrak{B}^C$ forms a basis of $\mathrm{DY}_{h}(\mathfrak{gl}^{\mathfrak{s}}_{m|n})$ over $\mathcal{A}$.
\end{theo}

\begin{rema}
Note that the total order defined on the generators of $\mathrm{DY}_{h}(\mathfrak{gl}^{\mathfrak{s}}_{m|n})$ differs from that introduced in \cite{BK25-2}. Consequently, when $\mathfrak{s}=\mathfrak{s}^{st}$, the explicit form of the PBW basis for $\mathrm{DY}_{h}(\mathfrak{gl}^{\mathfrak{s}}_{m|n})$ constructed here deviates from the counterpart in \cite{BK25-2}. Nevertheless, the PBW basis from \cite{BK25-2} can be recovered by invoking the defining relations \eqref{Y1}–\eqref{Y3}.
\end{rema}

\subsection{Classical limit}
Now, we define a natural ascending filtration on $\mathrm{DY}_{h}(\mathfrak{gl}^{\mathfrak{s}}_{m|n})$ by setting
$\operatorname{deg}\,t_{ij}^{(r)}=r-1,$ $\operatorname{deg}\,t_{ij}^{(-r)}=-r$ for all $r\in\mathbb{Z}_{+}$, and $\operatorname{deg}\,C=\operatorname{deg}\,h=0$.
Denote by $\bar{t}_{ij}^{(\pm r)}$ the image of $t_{ij}^{(\pm r)}$ in the corresponding graded component of the associated graded superalgebra
$\operatorname{gr} \mathrm{DY}_{h}(\mathfrak{gl}^{\mathfrak{s}}_{m|n})$.
The following result holds for the graded superalgebra.

\begin{prop} \label{Iso 1}
The assignment
$$
E_{ij}^{(r-1)}\mapsto d_{i}\bar{t}_{ij}^{(r)},\quad E_{ij}^{(-r)}\mapsto d_{i}\bar{t}_{ij}^{(-r)}, \quad K\mapsto \bar{C},\quad h\mapsto \bar{h},
$$
defines an $\mathcal{A}$-superalgebra isomorphism
\begin{align}\label{map1}
\mathrm{U}(\widehat{\mathfrak{gl}}^{\mathfrak{s}}_{m|n})[[h]]\rightarrow \operatorname{gr} \mathrm{DY}_{h}(\mathfrak{gl}^{\mathfrak{s}}_{m|n}).
\end{align}
\end{prop}
\begin{proof}
From \eqref{Y1}-\eqref{Y3}, we obtain the following commutation relations in $\operatorname{gr}\mathrm{DY}_{h}(\mathfrak{gl}^{\mathfrak{s}}_{m|n}):$
\begin{align*}
&\left[d_{i}\bar{t}_{ij}^{(r)},d_{k}\bar{t}_{kl}^{(s)}\right]=\delta_{kj}d_{i}
\bar{t}_{il}^{(r+s-1)}-\delta_{il}\epsilon_{ij;kl}d_{k}\bar{t}_{kj}^{(r+s-1)},\\
&\left[d_{i}\bar{t}_{ij}^{(-r)},d_{k}\bar{t}_{kl}^{(-s)}\right]=\delta_{kj}d_{i}
\bar{t}_{il}^{(-r-s)}-\delta_{il}\epsilon_{ij;kl}d_{k}\bar{t}_{kj}
^{(-r-s)},\\
&\left[d_{i}\bar{t}_{ij}^{(r)},d_{k}\bar{t}_{kl}^{(-s)}\right]=
\left\{
  \begin{array}{ll}
   \delta_{kj}d_{i}
\bar{t}_{il}^{(r-s-1)}-\delta_{il}\epsilon_{ij;kl}d_{k}\bar{t}_{kj}
^{(r-s-1)}+(r-1)d_{i}\delta_{kj}\delta_{il}\delta_{r,s+1}\bar{C}, & r\leq s; \\
~\\
    \delta_{kj}d_{i}
\bar{t}_{il}^{(r-s)}-\delta_{il}\epsilon_{ij;kl}d_{k}\bar{t}_{kj}
^{(r-s)}+(r-1)d_{i}\delta_{kj}\delta_{il}\delta_{r,s+1}\bar{C}, & r>s.
  \end{array}
\right.
\end{align*}
According to \eqref{affine}, these relations imply that the map \eqref{map1} is a surjective superalgebra homomorphism. Furthermore, the injectivity of this map follows from Theorem \ref{PBW basis} and the PBW theorem for $U(\widehat{\mathfrak{gl}}_{m|n}^{\mathfrak{s}})$.
\end{proof}
\begin{rema}
Proposition \ref{Iso 1} implies that $\mathrm{DY}_{h}(\mathfrak{gl}_{m|n}^{\mathfrak{s}})$ is a flat deformation of $\mathrm{U}(\widehat{\mathfrak{gl}}_{m|n}^{\mathfrak{s}})$.
\end{rema}

\medskip
\section{Gaussian generators for the double super Yangians}

In this section, for any parity $\mathfrak{s}\in \mathcal{S}_{m|n}$, we describe a decomposition of the generating matrices $T^{\pm}(u)$ of the double super Yangian $\mathrm{DY}_{h}(\mathfrak{gl}_{m|n}^{\mathfrak{s}})$ via the quasideterminant approach \cite{GGRW05}. This decomposition provides the Gaussian generators for $\mathrm{DY}_{h}(\mathfrak{gl}_{m|n}^{\mathfrak{s}})$, and we further deduce the commutation relations among them. The results in this section will be used in the next section to define a new presentation for $\mathrm{DY}_{h}(\mathfrak{gl}_{m|n}^{\mathfrak{s}})$.

\subsection{Quasideterminants and Gauss decomposition}

Let $X=(x_{ij})$ be an $N\times N$ matrix  over a unital ring. For $i,j\in I_{\mathfrak{s}}$, let $X^{ij}$ be the submatrix obtained by removing the $i$-th row and $j$-th column of $X$, and assume that $X^{ij}$ is invertible.
The $(i,j)$-th quasideterminant of $X$ is defined by
$$
|X|_{i j}=x_{ij}-r_{i}^{j}(X^{ij})^{-1}c_{j}^{i}.
$$
Here, $r_{i}^{j}$ denotes the row vector from the $i$-th row of $X$ with $x_{ij}$ deleted, and $c_{i}^{j}$ denotes the column vector from the $j$-th column of $X$ with $x_{ij}$ deleted; for more details, see \cite{GGRW05}.
The quasideterminant $|X|_{i j}$ also has a graphical notation as below:
$$
|X|_{i j}=\left|\begin{array}{ccccc}
x_{11} & \cdots & x_{1 j} & \cdots & x_{1 N} \\
\cdots & & \cdots & \\
x_{i 1} & \cdots & \boxed{x_{i j}} & \cdots & x_{i N} \\
\cdots & & \cdots & \\
x_{N 1} & \cdots & x_{N j} & \cdots & x_{N N}
\end{array}\right|.
$$

Now we introduce the \textit{Gaussian generators} for the superalgebra $\mathrm{DY}_{h}(\mathfrak{gl}^{\mathfrak{s}}_{m|n})$ by the following quasideterminantal formulae. For $i,j\in I_{\mathfrak{s}}$ with $i<j$, define
\begin{align*}
k_i^{\pm}(u)&=\left|\begin{array}{cccc}
t_{11}^{\pm}(u) & \cdots & t_{1, i-1}^{\pm}(u) & t_{1 i}^{\pm}(u) \\
\vdots & \ddots & \vdots& \vdots \\
t_{i 1}^{\pm}(u) & \cdots & t_{i, i-1}^{\pm}(u) & \boxed{t_{i i}^{\pm}(u)}
\end{array}\right|,
\end{align*}
\begin{align*}
e_{i j}^{\pm}(u)&={k_i^{\pm}(u)}^{-1}\left|\begin{array}{cccc}
t_{11}^{\pm}(u) & \cdots & t_{1, i-1}^{\pm}(u) & t_{1 j}^{\pm}(u) \\
\vdots & \ddots & \vdots & \vdots \\
t_{i-1, i}^{\pm}(u) & \cdots & t_{i-1, i-1}^{\pm}(u) & t_{i-1, j}^{\pm}(u) \\
t_{i 1}^{\pm}(u) & \cdots & t_{i, i-1}^{\pm}(u) & \boxed{t_{i j}^{\pm}(u)}
\end{array}\right|,
\end{align*}
\begin{align*}
f_{j i}^{\pm}(u)&=\left|\begin{array}{cccc}
t_{11}^{\pm}(u) & \cdots & t_{1, i-1}^{\pm}(u) & t_{1 i}^{\pm}(u) \\
\vdots & \ddots & \vdots & \vdots \\
t_{i-1,1}^{\pm}(u) & \cdots & t_{i-1, i-1}^{\pm}(u) & t_{i-1, i}^{\pm}(u) \\
t_{j i}^{\pm}(u) & \cdots & t_{j, i-1}^{\pm}(u) & \boxed{t_{j i}^{\pm}(u)}
\end{array}\right|{k_i^{\pm}(u)}^{-1}.
\end{align*}

It is straightforward to check that the coefficients of the above series have the following expansions:
\begin{align*}
&f_{ji}^{+}(u)=h\sum_{r\geqslant 1}f_{ji}^{(r)}u^{-r},\quad e_{ij}^{+}(u)=h\sum_{r\geqslant 1}e_{ij}^{(r)}u^{-r},\quad k_{i}^{+}(u)=1+h\sum_{r\geqslant 1}k_{i}^{(r)}u^{-r},\\
&f_{ji}^{-}(u)=-h\sum_{r\geqslant 1}f_{ji}^{(-r)}u^{r-1},\quad e_{ij}^{-}(u)=-h\sum_{r\geqslant 1}e_{ij}^{(-r)}u^{r-1},\quad k_{i}^{-}(u)=1-h\sum_{r\geqslant 1}k_{i}^{(-r)}u^{r-1},
\end{align*}
where all coefficients $f_{ji}^{(\pm r)}$, $e_{ij}^{(\pm r)}$, and $k_{i}^{(\pm r)}$ lie in $\mathrm{DY}_{h}(\mathfrak{gl}_{m|n}^{\mathfrak{s}})$. In the following, we use the abbreviated notation for $i\in I_{\mathfrak{s}}'$:
\begin{align*}
    &e_{i}^{+}(u)=h\sum_{r\geqslant 1}e_{i}^{(r)}u^{-r}:=e_{i,i+1}^{+}(u),~ f_{i}^{+}(u)=h\sum_{r\geqslant 1}f_{i}^{(r)}u^{-r}:=f_{i+1,i}^{+}(u)\in \mathrm{DY}_{h}(\mathfrak{gl}_{m|n}^{\mathfrak{s}})[[u^{-1}]], \\
    &e_{i}^{-}(u)=-h\sum_{r\geqslant 1}e_{i}^{(-r)}u^{r-1}:=e_{i,i+1}^{-}(u),~f_{i}^{-}(u)=-h\sum_{r\geqslant 1}f_{i}^{(-r)}u^{r-1}:=f_{i+1,i}^{-}(u)\in \mathrm{DY}_{h}(\mathfrak{gl}_{m|n}^{\mathfrak{s}})[[u]].
\end{align*}
For $1<i+1<j\leqslant N$, we have the recursive relations:
\begin{align}\label{ef}
e_{ij}^{(\pm r)}=d_{j-1}\left[e_{i,j-1}^{(\pm r)},e_{j-1}^{(1)}\right],\quad~
f_{ji}^{(\pm r)}=d_{j-1}\left[f_{j-1}^{(1)},f_{j-1,i}^{(\pm r)}\right].
\end{align}

By definition, the leading principal minors of the matrix $T^{\pm}(u)$ are invertible. Thus, $T^{\pm}(u)$ admit the following Gauss decomposition.

\begin{prop}\label{Gauss decomposition}
The generating matrices $T^{\pm}(u)$ have the unique Gauss decomposition:
\begin{equation}\label{Gauss_dec}
    \begin{split}
        T^{\pm}(u)&=\left(\begin{array}{cccc}
        1 & & & 0 \\
        f_{21}^{\pm}(u) & \ddots & &  \\
        \vdots & & \ddots \\
        f_{N, 1}^{\pm}(u) & f_{N, 2}^{\pm}(u) & \cdots & 1
        \end{array}\right)
        \left(\begin{array}{cccc}
        k_1^{\pm}(u) & &  & 0 \\
        & k_2^{\pm}(u) & &  \\
         & & \ddots & \\
        0 &  & & k^{\pm}_{N}(u)
        \end{array}\right)\\ 
        &\hspace{3em}\times\left(\begin{array}{cccc}
        1 & e_{12}^{\pm}(u) & \cdots & e_{1, N}^{\pm}(u) \\
        & \ddots & & e_{2, N}^{\pm}(u) \\
        & & \ddots & \vdots \\
        0 & & & 1
        \end{array}\right),
    \end{split}
\end{equation}
where $k_{i}^{\pm}(u)$ are invertible. Thus $\mathrm{DY}_{h}(\mathfrak{gl}^{\mathfrak{s}}_{m|n})$ is generated by the coefficients of series
\begin{align*}
\{k_{i}^{\pm}(u), e_{j}^{\pm}(u), f_{j}^{\pm}(u)\mid i\in I_{\mathfrak{s}}; j\in I_{\mathfrak{s}}'\}.
\end{align*}
\end{prop}

More precisely, the entries of $T^{\pm}(u)$ satisfy the following relations for  $i,j\in I_{\mathfrak{s}}$ with $i<j$:
\begin{align}\label{Gauss ii}
t_{ii}^{\pm}(u)&= k_{i}^{\pm}(u) + \sum_{s < i} f_{is}^{\pm}(u) k_{s}^{\pm}(u) e_{si}^{\pm}(u),\\ \label{Gauss ij}
t_{ij}^{\pm}(u)&= k_{i}^{\pm}(u) e_{ij}^{\pm}(u) + \sum_{s< i} f_{is}^{\pm}(u) k_{s}^{\pm}(u) e_{sj}^{\pm}(u),\\ \label{Gauss ji}
t_{ji}^{\pm}(u)&= f_{ji}^{\pm}(u) k_{i}^{\pm}(u) + \sum_{s <i} f_{js}^{\pm}(u) k_{s}^{\pm}(u) e_{si}^{\pm}(u).
\end{align}

\subsection{Gaussian relations}

To explicitly write down the relations among the Gauss generators, we start with the special cases where $N=m+n$ is either 2 or 3, which are less complicated. The other relations in full generality can be deduced from these special ones by applying certain injective homomorphisms introduced in the previous subsection.

Set for $i\in I_{\mathfrak{s}}'$,
\begin{equation*}
X_{i}^{-}(u)=f_{i}^{+}(u_{+})-f_{i}^{-}(u_{-}),\qquad X_{i}^{+}(u)=e_{i}^{+}(u_{-})-e_{i}^{-}(u_{+}).
\end{equation*}
 Notice that $|X_i^\pm(u)| = |\alpha_i| \equiv s_i+s_{i+1} \pmod 2$.

\subsubsection{The case of $N=2$}
By definition, we have the following formulas:
\begin{align*}
&T^{ \pm}(u)=\left(\begin{array}{cc}
k_1^{ \pm}(u) & k_1^{ \pm}(u) e_1^{ \pm}(u) \\
f_1^{ \pm}(u) k_1^{ \pm}(u) & k_2^{ \pm}(u)+f_1^{ \pm}(u) k_1^{ \pm}(u) e_1^{ \pm}(u)
\end{array}\right),\\
&T^{\pm}(u)^{-1}=\left(\begin{array}{cc}
k_1^{\pm}(u)^{-1}+e_1^{\pm}(u) k_2^{\pm}(u)^{-1}f_1^{ \pm}(u) & -e_1^{\pm}(u)k_2^{ \pm}(u)^{-1} \\
-k_2^{\pm}(u)^{-1}f_1^{\pm}(u) & k_2^{\pm}(u)^{-1}
\end{array}\right),\\
&R_{21}(u)=R_{12}(u)=\left( {\begin{array}{*{20}{c}}
{1+\frac{d_{1}}{u}h} & 0 & 0 & 0 \\
0 & 1 & {\frac{d_{2}}{u}h} & 0 \\
0 & {\frac{d_{1}}{u}h} & 1 & 0 \\
0 & 0 & 0 & {1+\frac{d_{2}}{u}h}
\end{array}}\right).
\end{align*}

\begin{rema}\label{m,n}
For $N=2$, we have three separate cases: $(m,n)=(2,0)$, $(m,n)=(1,1)$ and $(m,n)=(0,2)$. The associated R-matrices only differ in the values of $d_{i},~i=1,2$. Hence all cases can be handled uniformly, with attention paid merely to the values of $d_{i}$ case by case.
\end{rema}
Following the same argument as in the proof of Theorem 2.5 in \cite{YJ24} for the double Yangian $\mathrm{DY}_{h}(\mathfrak{gl}_{n})$, the following relations are obtained from \eqref{RTT1}, \eqref{RTT2} and \eqref{RTT4}-\eqref{RTT9}:
{\allowdisplaybreaks
\begin{align}
k_{i}^{\pm}(u)k_{i}^{\pm}(v)&=k_{i}^{\pm}(v)k_{i}^{\pm}(u),\quad i=1,2,\\
k_{2}^{\pm}(v)k_{1}^{\pm}(u)&=k_{1}^{\pm}(u)k_{2}^{\pm}(v),\\
k_{2}^{-}(v)k_{1}^{+}(u)&=k_{1}^{+}(u)k_{2}^{-}(v),\\
\left(1+\frac{d_{i}}{u_{-}-v_{+}}h\right)k_{i}^{+}(u)k_{i}^{-}(v)&=
\left(1+\frac{d_{i}}{u_{+}-v_{-}}h\right)k_{i}^{-}(v)k_{i}^{+}(u),~i=1,2,\\
\frac{(u_{-}-v_{+})^{2}}{(u_{-}-v_{+})^{2}-h^{2}}k_{1}^{-}(u)k_{2}^{+}(v)^{-1}&=
\frac{(u_{+}-v_{-})^{2}}{(u_{+}-v_{-})^{2}-h^{2}}k_{2}^{+}(v)^{-1}k_{1}^{-}(u),\\
k_{1}^{\pm}(u)^{-1}X_{1}^{+}(v)k_{1}^{\pm}(u)&=\left(1+\frac{d_{1}}{u_{\pm}-v}h\right) X_{1}^{+}(v),\\
k_{2}^{\pm}(u)^{-1}X_{1}^{+}(v)k_{2}^{\pm}(u)&=\left(1-\frac{d_{2}}{u_{\pm}-v}h\right) X_{1}^{+}(v),\\
k_{1}^{\pm}(u)X_{1}^{-}(v)k_{1}^{\pm}(u)^{-1}&=\left(1+\frac{d_{1}}{u_{\mp}-v}h \right)X_{1}^{-}(v),\\
k_{2}^{\pm}(u)X_{1}^{-}(v)k_{2}^{\pm}(u)^{-1}&=\left(1-\frac{d_{2}}{u_{\mp}-v}h \right)X_{1}^{-}(v),\\ \label{X1X1+}
(-1)^{|\alpha_{1}|}\left(u-v-d_{1}h\right)X_{1}^{+}(u)X_{1}^{+}(v)&=
\left(u-v+d_{2}h\right)X_{1}^{+}(v)X_{1}^{+}(u),\\ \label{X1X1-}
(-1)^{|\alpha_{1}|}\left(u-v-d_{1}h\right)X_{1}^{-}(v)X_{1}^{-}(u)&=
\left(u-v+d_{2}h\right)X_{1}^{-}(u)X_{1}^{-}(v),\\
\left[X_{1}^{+}(u),X_{1}^{-}(v)\right]=(-1)^{|1||2|}h\{\delta(\frac{v_{+}}{u_{-}})&k_{2}^{+}(u_{-})k_{1}^{+}(u_{-})^{-1}-\delta(\frac{v_{-}}{u_{+}})k_{2}^{-}(v_{-})k_{1}^{-}(v_{-})^{-1}\},
\end{align}
}
where $\delta(\frac{v}{u})=\sum_{k\in \mathbb{Z}}u^{-k-1}v^{k}$.

\subsubsection{The case of $N=3$}
For $N=3$, there exist four cases, all admitting uniform treatment as noted in Remark \ref{m,n}. In this subsection, we write $R_{N}(u)$ for the R-matrix $R(u)$ associated with dimension $N$. This notation should not be confused with $R_{21}(u)$ and similar symbols used in previous sections.

Restricting the relations \eqref{RTT1} and \eqref{RTT2} to $E_{ij}\otimes E_{kl}$, with $1\leqslant i,j,k,l\leqslant 2$, we obtain
\begin{align*}
R_{2}(u-v)L_{1}^{\pm}(u)L_{2}^{\pm}(v)&=L_{2}^{\pm}(v)L_{1}^{\pm}(u)R_{2}(u-v),\\ 
R_{2}(u_{-}-v_{+})L_{1}^{+}(u)L_{2}^{-}(v)&=L_{2}^{-}(v)L_{1}^{+}(u)R_{2}(u_{+}-v_{-}),
\end{align*}
where 
\begin{align*}
L^{\pm}(u)=\left(\begin{array}{cc}
1 & 0 
\\
f_1^{\pm}(u) & 1
\end{array}\right)
\left(\begin{array}{cc}
k_{1}^{\pm}(u)& 0 
\\
0 & k_{2}^{\pm}(u)
\end{array}\right)
\left(\begin{array}{cc}
1 & e_1^{\pm}(u)
\\
0 & 1
\end{array}\right).
\end{align*}
These are precisely the relations of the $N=2$ case. Similarly, restricting \eqref{RTT1} and \eqref{RTT2} to $E_{ij}\otimes E_{kl}$, with $2\leqslant i,j,k,l\leqslant 3$ yields
\begin{align*}
R_{2}(u-v)\widehat{L}_{1}^{\pm}(u)\widehat{L}_{2}^{\pm}(v)&=\widehat{L}_{2}^{\pm}(v)\widehat{L}_{1}^{\pm}(u)R_{2}(u-v),\\ 
R_{2}(u_{-}-v_{+})\widehat{L}_{1}^{+}(u)\widehat{L}_{2}^{-}(v)&=\widehat{L}_{2}^{-}(v)\widehat{L}_{1}^{+}(u)R_{2}(u_{+}-v_{-}),
\end{align*}
where
\begin{align*}
\widehat{L}^{\pm}(u)=\left(\begin{array}{cc}
1 & 0 
\\
f_2^{\pm}(u) & 1
\end{array}\right)
\left(\begin{array}{cc}
k_{2}^{\pm}(u)& 0 
\\
0 & k_{3}^{\pm}(u)
\end{array}\right)
\left(\begin{array}{cc}
1 & e_2^{\pm}(u)
\\
0 & 1
\end{array}\right). 
\end{align*}
This is analogous to the $N=2$ case. Consequently, it remains to derive the relations between the sets $\{k_1^{\pm}(u)$, $e_1^{\pm}(u)$, $f_1^{\pm}(u)\}$ and $\{k_3^{\pm}(u), e_2^{\pm}(u)$, $f_2^{\pm}(u)\}$.
By the Gauss decomposition of $T^{\pm}(u)$, we have
$$
T^{\pm}(u)=\left(\begin{array}{ccc}
k_1^{\pm}(u) & k_1^{\pm}(u)e_1^{\pm}(u) & k_{1}^{\pm}(u)e_{1,3}^{\pm}(u)
\\
f_1^{\pm}(u)k_1^{\pm}(u) & * & *
\\
f_{3,1}^{\pm}(u) k_1^{\pm}(u) & * & *
\end{array}\right).
$$
Let $x^{\pm}=k_{3}^{\pm}(v)^{-1}(-f_{3,1}^{\pm}(v)+f_{2}^{\pm}(v)f_{1}^{\pm}(v))$ and
 $y^{\pm}=(-e_{1,3}^{\pm}(v)+e_{1}^{\pm}(v)e_{2}^{\pm}(v))k_{3}^{\pm}(v)^{-1}$, then
$$
T^{\pm}(v)^{-1}=\left(\begin{array}{ccc}
* & * & y^{\pm}
\\
* & * & -e_2^{\pm}(v)k_3^{\pm}(v)^{-1}
\\
x^{\pm} & -k_3^{\pm}(v)^{-1}f_2^{\pm}(v) & k_3^{\pm}(v)^{-1}
\end{array}\right),
$$
where $*$ represent some elements in $\mathrm{DY}_{h}(\mathfrak{gl}^{\mathfrak{s}}_{m|n})$.

\begin{lemm}
We have the following relations in $\mathrm{DY}_{h}(\mathfrak{gl}^{\mathfrak{s}}_{m|n}):$
\begin{align}\label{k3k1}
k_{3}^{\pm}(v)^{-1}k_{1}^{\pm}(u)&=
k_{1}^{\pm}(u)k_{3}^{\pm}(v)^{-1},\\ \label{k3-k1+}
k_{3}^{-}(v)^{-1}k_{1}^{+}(u)&=k_{1}^{+}(u)k_{3}^{-}(v)^{-1},\\ \label{k3+k1-}
\frac{(u_{+}-v_{-})^{2}}{(u_{+}-v_{-})^{2}-h^{2}}k_{3}^{+}(v)^{-1}k_{1}^{-}(u)&=
\frac{(u_{-}-v_{+})^{2}}{(u_{-}-v_{+})^{2}-h^{2}}k_{1}^{-}(u)k_{3}^{+}(v)^{-1},\\ \label{k1x2}
k_{1}^{\iota}(u)^{-1}
X_{2}^{\iota'}(v)k_{1}^{\iota}(u)&=X_{2}^{\iota'}(v),\\ \label{k3x1}
k_{3}^{\iota}(u)^{-1}
X_{1}^{\iota'}(v)k_{3}^{\iota}(u)&=X_{1}^{\iota'}(v),\\
\label{x1x2}
\left[X_{2}^{+}(u),X_{1}^{-}(v)\right]&=\left[X_{1}^{+}(u),X_{2}^{-}(v)\right]=0,\\ \label{X1X2-}
X_{1}^{-}(u)X_{2}^{-}(v)=(-1)^{|\alpha_{1}||\alpha_{2}|}&\left(1+\frac{d_{2}}{u-v}h\right)X_{2}^{-}(v)X_{1}^{-}(u),\\ \label{X2X1+}
X_{2}^{+}(v)X_{1}^{+}(u)=(-1)^{|\alpha_{1}||\alpha_{2}|}&\left(1+\frac{d_{2}}{u-v}h\right)X_{1}^{+}(u)X_{2}^{+}(v),\\ \notag
X_{1}^{\pm}(u_{1})X_{1}^{\pm}(u_{2})X_{2}^{\pm}(v)-2
X_{1}^{\pm}(u_{1})&X_{2}^{\pm}(v)X_{1}^{\pm}(u_{2})+X_{2}^{\pm}(v)X_{1}^{\pm}(u_{1})X_{1}^{\pm}(u_{2})\\
\label{serre3-1}
&+\{u_{1}\leftrightarrow u_{2}\}=0, \quad~|\alpha_{1}|=0,\\ \notag
X_{2}^{\pm}(u_{1})X_{2}^{\pm}(u_{2})X_{1}^{\pm}(v)-2 X_{2}^{\pm}(u_{1})&X_{1}^{\pm}(v)X_{2}^{\pm}(u_{2})+X_{1}^{\pm}(v)X_{2}^{\pm}(u_{1})X_{2}^{\pm}(u_{2})\\ \label{serre3-2}
&+\{u_{1}\leftrightarrow u_{2}\}=0, \quad~|\alpha_{2}|=0,
\end{align}
where $\iota,\iota'\in\{+,-\}$.
\end{lemm}
\begin{proof}
From \eqref{RTT4}, \eqref{RTT5} and \eqref{RTT9}, we can obtain \eqref{k3k1}-\eqref{k3+k1-} and
\begin{align}\label{k3e1}
k_{3}^{\iota}(v)^{-1}e_1^{\iota'}(u)&=e_1^{\iota'}(u)k_{3}^{\iota}(v)^{-1},\\ \label{k3f1}
k_{3}^{\iota}(v)^{-1}f_{1}^{\iota'}(u)&=f_{1}^{\iota'}(u)k_{3}^{\iota}(v)^{-1},\\ \label{k1e2}
e_{2}^{\iota}(v)k_{1}^{\iota'}(u)&=k_{1}^{\iota'}(u)e_{2}^{\iota}(v),\\ \label{k1f2}
f_{2}^{\iota}(v)k_{1}^{\iota'}(u)&=k_{1}^{\iota'}(u)f_{2}^{\iota}(v),\\ \label{e2f1+}
e_{2}^{\iota}(v)f_{1}^{\iota'}(u)&=(-1)^{|\alpha_{1}||\alpha_{2}|}
f_{1}^{\iota'}(u)e_{2}^{\iota}(v),\\  \label{f2e1-}
f_{2}^{\iota}(v)e_{1}^{\iota'}(u)&=(-1)^{|\alpha_{1}||\alpha_{2}|}
e_{1}^{\iota'}(u)f_{2}^{\iota}(v).
\end{align}
By \eqref{k3e1}-\eqref{f2e1-} together with the definition of $X_{i}^{\pm}, i=1,2$, we can deduce \eqref{k1x2}-\eqref{x1x2}.

We now verify the relations among $X_{1}^{\pm}(u)$ and $X_{2}^{\mp}(v)$. From \eqref{RTT4}, \eqref{RTT5} and \eqref{RTT9}, we have
\begin{align*}
f_{1}^{\pm}(u)f_{2}^{\pm}(v)=&-(-1)^{|1||3|+|2||3|}\frac{d_{2}h}{u-v}
\left((-1)^{|1||2|+|2||3|}f_{3,1}^{\pm}(u)-f_{3,1}^{\pm}(v)\right)\\
-&(-1)^{|1||3|+|2||3|}\frac{d_{2}h}{u-v}
f_{2}^{\pm}(v)f_{1}^{\pm}(v)+(-1)^{|\alpha_{1}||\alpha_{2}|}\left(1+\frac{d_{2}h}{u-v}\right)
f_{2}^{\pm}(v)f_{1}^{\pm}(u),\\
f_{1}^{\mp}(u)f_{2}^{\pm}(v)=&-(-1)^{|1||3|+|2||3|}\frac{d_{2}h}{u_{\pm}-v_{\mp}}
\left((-1)^{|1||2|+|2||3|}f_{3,1}^{\mp}(u)-f_{3,1}^{\pm}(v)\right)\\
-&(-1)^{|1||3|+|2||3|}\frac{d_{2}h}
{u_{\pm}-v_{\mp}}f_{2}^{\pm}(v)f_{1}^{\pm}(v)+(-1)^{|\alpha_{1}||\alpha_{2}|}\left(1+\frac{d_{2}h}{u_{\pm}-v_{\mp}}\right)
f_{2}^{\pm}(v)f_{1}^{\mp}(u).
\end{align*}
By the definition of $X_{1}^{-}(u)$ and $X_{2}^{-}(v)$, we get \eqref{X1X2-}. Similarly, we can prove that \eqref{X2X1+} holds.

Finally, we check the Serre relations \eqref{serre3-1} and \eqref{serre3-2}. Indeed, the verifications of \eqref{serre3-1} and \eqref{serre3-2} are similar; we only consider the case of \eqref{serre3-1} as an example. When $|\alpha_{1}|=0$, from \eqref{X1X1-} and \eqref{X1X2-}, we get
\begin{align*}
&X_{1}^{-}(u_{1})X_{2}^{-}(v)X_{1}^{-}(u_{2})=
\frac{u_{2}-v}{u_{2}-v+d_{2}h}X_{1}^{-}(u_{1})X_{1}^{-}(u_{2})X_{2}^{-}(v),\\
&X_{2}^{-}(v)X_{1}^{-}(u_{1})X_{1}^{-}(u_{2})=\frac{u_{2}-v}{u_{2}-v+d_{2}h}
\frac{u_{1}-v}{u_{1}-v+d_{2}h}X_{1}^{-}(u_{1})X_{1}^{-}(u_{2})X_{2}^{-}(v),\\
&X_{1}^{-}(u_{2})X_{1}^{-}(u_{1})X_{2}^{-}(v)=\frac{u_{1}-u_{2}+d_{2}h}
{u_{1}-u_{2}-d_{1}h}X_{1}^{-}(u_{1})X_{1}^{-}(u_{2})X_{2}^{-}(v),\\
&X_{1}^{-}(u_{2})X_{2}^{-}(v)X_{1}^{-}(u_{1})=\frac{u_{1}-v}{u_{1}-v+d_{2}h}\frac{u_{1}-u_{2}+d_{2}h}{u_{1}-u_{2}-d_{1}h}
X_{1}^{-}(u_{1})X_{1}^{-}(u_{2})X_{2}^{-}(v),\\
&X_{2}^{-}(v)X_{1}^{-}(u_{2})X_{1}^{-}(u_{1})=\frac{u_{1}-v}
{u_{1}-v+d_{2}h}\frac{u_{2}-v}{u_{2}-v+d_{2}h}\frac{u_{1}-u_{2}+d_{2}h}
{u_{1}-u_{2}-d_{1}h}X_{1}^{-}(u_{1})X_{1}^{-}(u_{2})X_{2}^{-}(v).
\end{align*}
Hence, each summand in \eqref{serre3-1} can be replaced by $X_{1}^{-}(u_{1})X_{1}^{-}(u_{2})X_{2}^{-}(v)$. A direct computation then shows that the coefficients sum to zero, so \eqref{serre3-1} holds.
\end{proof}

\subsubsection{General cases}

Now we proceed to the general $N\geqslant4$. Arguing by induction, we assume that all relations for the case $N-1$ are established. As in the case $N=3$, it is sufficient to verify the relations between the sets $\{k_1^{\pm}(u), e_1^{\pm}(u), f_1^{\pm}(u)\}$ and $\{k_N^{\pm}(u)$, $e_{N-1}^{\pm}(u), f_{N-1}^{\pm}(u)\}$.

\begin{lemm}
The following relations hold in $\mathrm{DY}_{h}(\mathfrak{gl}^{\mathfrak{s}}_{m|n}):$
\begin{align}\label{kNk1}
&k_{N}^{\pm}(v)^{-1}k_{1}^{\pm}(u)=k_{1}^{\pm}(u)k_{N}^{\pm}(v)^{-1},\\ \label{kN-k1+}
&k_{N}^{-}(v)^{-1}k_{1}^{+}(u)=k_{1}^{+}(u)k_{N}^{-}(v)^{-1},\\ \label{kN+k1-}
&\frac{(u_{+}-v_{-})^{2}}{(u_{+}-v_{-})^{2}-h^{2}}k_{N}^{+}(v)^{-1}k_{1}^{-}(u)=\frac{(u_{-}-v_{+})^{2}}{(u_{-}-v_{+})^{2}-h^{2}}k_{1}^{-}(u)k_{N}^{+}(v)^{-1},\\ \label{k1XN-1}
&k_{1}^{\iota'}(u)^{-1}X_{N-1}^{\iota}(v)k_{1}^{\iota'}(u)=X_{N-1}^{\iota}(v),\\ \label{kNX1}
&k_{N}^{\iota'}(u)^{-1}X_{1}^{\iota}(v)k_{N}^{\iota'}(u)=X_{1}^{\iota}(v),\\
\label{X1XN-1}
&\left[X_{1}^{\iota'}(u),X_{N-1}^{\iota}(v)\right]=0,\\ \label{serre i}
&{\rm Sym}_{u_{1},u_{2}}\left[X_{i}^{\pm}(u_{1}),\left[X_{i+1}^{\pm}(v_{1}),\left[
X_{i}^{\pm}(u_{2}),X_{i-1}^{\pm}(v_{2})\right]\right]\right]=0, \quad |\alpha_{i}|=1, i\pm1\in I_{\mathfrak{s}}',
\end{align}
where $\iota,\iota'\in\{+,-\}$.
\end{lemm}
\begin{proof}
By \eqref{RTT4}, \eqref{RTT5} and \eqref{RTT9}, we obtain \eqref{kNk1}-\eqref{kN+k1-} and the following relations:
\begin{align}
\label{kNe1}
k_{N}^{\iota}(v)^{-1}e_{1}^{\iota'}(u)&=e_{1}^{\iota'}(u)k_{N}^{\iota}(v)^{-1},\\  \label{kNf1}
k_{N}^{\iota}(v)^{-1}f_{1}^{\iota'}(u)&=f_{1}^{\iota'}(u)k_{N}^{\iota}(v)^{-1},\\ \label{k1fN-1}
f_{N-1}^{\iota}(v)k_{1}^{\iota'}(u)&=k_{1}^{\iota'}(u)f_{N-1}^{\iota}(v),\\ \label{k1eN-1}
e_{N-1}^{\iota}(v)k_{1}^{\iota'}(u)
&=k_{1}^{\iota'}(u)e_{N-1}^{\iota}(v),\\ \label{eN-1f1}
e_{N-1}^{\iota}(v)
f_{1}^{\iota'}(u)&=(-1)^{|\alpha_{1}||\alpha_{N-1}|}
f_{1}^{\iota'}(u)e_{N-1}^{\iota}(v),\\  \label{e1fN-1}
f_{N-1}^{\iota}(v)e_{1}^{\iota'}(u)&=(-1)^{|\alpha_{1}||\alpha_{N-1}|}
e_{1}^{\iota'}(u)f_{N-1}^{\iota}(v),\\ \label{f1fN-1}
f_{N-1}^{\iota}(v)f_{1}^{\iota'}(u)&=(-1)^{|\alpha_{1}||\alpha_{N-1}|}
f_{1}^{\iota'}(u)f_{N-1}^{\iota}(v),\\ \label{e1eN-1}
e_{N-1}^{\iota}(v)
e_{1}^{\iota'}(u)&=(-1)^{|\alpha_{1}||\alpha_{N-1}|}
e_{1}^{\iota'}(u)e_{N-1}^{\iota}(v).
\end{align}
From \eqref{kNe1}-\eqref{e1eN-1}, we get \eqref{k1XN-1}-\eqref{X1XN-1}.

When $|\alpha_{i}|=1$, according to \eqref{X1X1+}, \eqref{X1X1-}, \eqref{X1X2-}, \eqref{X2X1+} and \eqref{X1XN-1}, for $i\pm1\in I_{\mathfrak{s}}$, we can deduce that
\begin{align*}
&X_{i}^{\pm}(u)X_{i}^{\pm}(v)=-X_{i}^{\pm}(v)X_{i}^{\pm}(u),\\
&X_{i+1}^{+}(v)X_{i}^{+}(u)=(-1)^{|\alpha_{i+1}||\alpha_{i}|}\left(1+\frac{d_{i+1}}{u-v}h\right)X_{i}^{+}(u)X_{i+1}^{+}(v),\\
&X_{i}^{-}(u)X_{i+1}^{-}(v)=(-1)^{|\alpha_{i+1}||\alpha_{i}|}\left(1+\frac{d_{i+1}}{u-v}h\right)X_{i+1}^{-}(v)X_{i}^{-}(u),\\
&\left[X_{i-1}^{\pm}(u),X_{i+1}^{\pm}(v)\right]=0.
\end{align*}
They yield the Serre relations:
\begin{align*}
&{\rm Sym}_{u_{1},u_{2}}\left[X_{i}^{\pm}(u_{1}),\left[X_{i+1}^{\pm}(v_{1}),\left[X_{i}^{\pm}(u_{2}), X_{i-1}^{\pm}(v_{2})\right]\right]\right]\\
=&X_{i}^{\pm}(u_{1})X_{i+1}^{\pm}(v_{1})X_{i}^{\pm}(u_{2})X_{i-1}^{\pm}(v_{2})+(-1)^{|\alpha_{i-1}|+|\alpha_{i+1}|}X_{i+1}^{\pm}(v_{1})X_{i}^{\pm}(u_{1})
X_{i-1}^{\pm}(v_{2})X_{i}^{\pm}(u_{2})\\
+&(-1)^{|\alpha_{i-1}|+|\alpha_{i+1}|+|\alpha_{i-1}||\alpha_{i+1}|}X_{i}^{\pm}(u_{1})
X_{i-1}^{\pm}(v_{2})X_{i}^{\pm}(u_{2})X_{i+1}^{\pm}(v_{1})\\
+&(-1)^{|\alpha_{i-1}|
|\alpha_{i+1}|}X_{i-1}^{\pm}(v_{2})X_{i}^{\pm}(u_{1})X_{i+1}^{\pm}(v_{1})
X_{i}^{\pm}(u_{2})\\
-&2\times(-1)^{|\alpha_{i-1}|+|\alpha_{i-1}||\alpha_{i+1}|}X_{i}^{\pm}(u_{1})
X_{i-1}^{\pm}(v_{2})X_{i+1}^{\pm}(v_{1})X_{i}^{\pm}(u_{2})\\
+&\{u_{1}\leftrightarrow u_{2}\}=0,
\end{align*}
whose verification is analogous to those of \eqref{serre3-1}–\eqref{serre3-2}.
\end{proof}

Next, we list all the relations among the currents $k_{i}^{\pm}(u)$, $X_{i}^{+}(u)$ and $X_{i}^{-}(v)$.

\begin{theo}\label{Drinfeld relation 1}
The following relations hold in the algebra $\mathrm{DY}_{h}(\mathfrak{gl}^{\mathfrak{s}}_{m|n}):$
\begin{align}\label{kikj}
&k_{i}^{\pm}(u)k_{j}^{\pm}(v)=k_{j}^{\pm}(v)k_{i}^{\pm}(u),\\ \label{ki+kj-}
&k_{i}^{+}(u)k_{j}^{-}(v)=k_{j}^{-}(v)k_{i}^{+}(u),\quad i<j,\\ \label{ki-kj+}
&\frac{(u_{+}-v_{-})^{2}}{(u_{+}-v_{-})^{2}-h^{2}}k_{i}^{-}(u)k_{j}^{+}(v)=
\frac{(u_{-}-v_{+})^{2}}{(u_{-}-v_{+})^{2}-h^{2}}k_{j}^{+}(v)k_{i}^{-}(u),\quad i<j,\\ \label{ki+ki-}
&\left(1+\frac{d_{i}}{u_{-}-v_{+}}h\right)k_{i}^{+}(u)k_{i}^{-}(v)=\left(1+
\frac{d_{i}}{u_{+}-v_{-}}h\right)k_{i}^{-}(v)k_{i}^{+}(u),\\ \label{kiXi+}
&k_{i}^{\pm}(u)^{-1}X_{i}^{+}(v)k_{i}^{\pm}(u)=\left(1+\frac{d_{i}}{u_{\pm}-v}h
\right)X_{i}^{+}(v),\\ \label{kiXi-}
&k_{i}^{\pm}(u)X_{i}^{-}(v)k_{i}^{\pm}(u)^{-1}=\left(1+\frac{d_{i}}{u_{\mp}-v}h
\right)X_{i}^{-}(v),\\ \label{ki+1Xi+}
&k_{i+1}^{\pm}(u)^{-1}X_{i}^{+}(v)k_{i+1}^{\pm}(u)=\left(1-\frac{d_{i+1}}
{u_{\pm}-v}h\right)X_{i}^{+}(v),\\ \label{ki+1Xi-}
&k_{i+1}^{\pm}(u)X_{i}^{-}(v)k_{i+1}^{\pm}(u)^{-1}=\left(1-\frac{d_{i+1}}
{u_{\mp}-v}h\right)X_{i}^{-}(v),\\ \label{kjXi}
&k_{j}^{\iota'}(u)^{-1}X_{i}^{\iota}(v)k_{j}^{\iota'}(u)=X_{i}^{\iota}(v), \quad j\neq i,i+1,\\ \label{XiXi}
&(-1)^{|\alpha_{i}|}\left(u-v\mp d_{i}h\right)X_{i}^{\pm}(u)X_{i}^{\pm}(v)=\left(u-v
\pm d_{i+1}h\right)X_{i}^{\pm}(v)X_{i}^{\pm}(u),\\ \label{XiXi+1-}
&X_{i}^{-}(u)X_{i+1}^{-}(v)=(-1)^{|\alpha_{i}||\alpha_{i+1}|}\left(1+\frac{d_{i+1}}
{u-v}h\right)X_{i+1}^{-}(v)X_{i}^{-}(u),\\ \label{XiXi+1+}
&X_{i+1}^{+}(v)X_{i}^{+}(u)=(-1)^{|\alpha_{i}||\alpha_{i+1}|}\left(1+\frac{d_{i+1}}
{u-v}h\right)X_{i}^{+}(u)X_{i+1}^{+}(v),\\ \label{XiXj}
&\left[X_{i}^{\pm}(u),X_{j}^{\pm}(v)\right]=0,\quad |i-j|>1,\\ \label{serre  i-1}
&{\rm Sym}_{u_{1},u_{2}}\left[X_{i}^{\pm}(u_{1}),\left[X_{i}^{\pm}(u_{2}),X_{i\pm 1}^{\pm}(v)\right]\right]=0, \quad |\alpha_{i}|=0,~i\pm1\in I_{\mathfrak{s}}',\\ \label{serre i-2}
&{\rm Sym}_{u_{1},u_{2}}\left[X_{i}^{\pm}(u_{1}),\left[X_{i+1}^{\pm}(v_{1}),\left[
X_{i}^{\pm}(u_{2}),X_{i-1}^{\pm}(v_{2})\right]\right]\right]=0, \quad |\alpha_{i}|=1,~i\pm1\in I_{\mathfrak{s}}',\\ \label{Xi+Xj-}
&\left[X_{i}^{+}(u),X_{j}^{-}(v)\right]=(-1)^{|i||i+1|}\delta_{ij}h\{\delta(\frac{v_{+}}{u_{-}})k_{i+1}^{+}(u_{-})k_{i}^{+}(u_{-})^{-1}-\delta(\frac{v_{-}}{u_{+}})k_{i+1}^{-}(v_{-})k_{i}^{-}(v_{-})^{-1}\},
\end{align}
where $\iota, \iota' \in\{+,-\}$ and $\delta(\frac{v}{u})=\sum_{k\in \mathbb{Z}}u^{-k-1}v^{k}$.
\end{theo}

\begin{rema}
In the case $\mathfrak{s}=\mathfrak{s}^{st}$, the relations of Theorem \ref{Drinfeld relation 1} are analogous to those appearing in \cite{Z97}, although the proof is not included there. We also note that the R-matrix $R(u)$ used here differs slightly from that in \cite{Z97}.
When $n=0$, these relations in Theorem \ref{Drinfeld relation 1} coincide with those in Theorem 2.5 of \cite{YJ24}.
\end{rema}

\medskip
\section{Drinfeld presentation of the double super Yangian}

In this section, we present the Drinfeld presentation of the double super Yangian associated with the Lie superalgebra $\mathfrak{gl}_{m|n}^{\mathfrak{s}}$ and establish an explicit isomorphism between the two presentations. Motivated by the idea of Yang and Jing \cite{YJ24}, we can prove the following theorem.

\begin{theo}\label{iso gl}
Let $\widehat{\mathrm{DY}}_{h}(\mathfrak{gl}_{m|n}^{\mathfrak{s}})$ be the $\mathcal{A}$-unital associative superalgebra  generated by the coefficients of the series $k_{i}^{\pm}(u), e_{j}^{\pm}(u)$ and $f_{j}^{\pm}(u)$ with $i\in I_{\mathfrak{s}}, j\in I_{\mathfrak{s}}'$ and the central element $C$. The defining relations are given by the relations in Theorem \ref{Drinfeld relation 1}. There exists a superalgebra isomorphism $\phi$ from $\widehat{\mathrm{DY}}_{h}(\mathfrak{gl}_{m|n}^{\mathfrak{s}})$ to  $\mathrm{DY}_{h}(\mathfrak{gl}_{m|n}^{\mathfrak{s}})$, which sends each generator of $\widehat{\mathrm{DY}}_{h}(\mathfrak{gl}_{m|n}^{\mathfrak{s}})$ to the corresponding generator of $\mathrm{DY}_{h}(\mathfrak{gl}_{m|n}^{\mathfrak{s}})$ with the same name.
\end{theo}
\begin{proof}
By Theorem \ref{Drinfeld relation 1}, the map $\phi$ is a superalgebra homomorphism. According to Proposition \ref{Gauss decomposition},  $\mathrm{DY}_{h}(\mathfrak{gl}^{\mathfrak{s}}_{m|n})$ is generated by the coefficients of the series $\{k_{i}^{\pm}(u), e_{j}^{\pm}(u), f_{j}^{\pm}(u)|~i\in I_{\mathfrak{s}}, j\in I_{\mathfrak{s}}'\}$ together with the central element $C$, so $\phi$ is surjective. It remains to prove injectivity. To this end, we show that $\widehat{\mathrm{DY}}_{h}(\mathfrak{gl}^{\mathfrak{s}}_{m|n})$ is spanned, as a vector space, by ordered monomials in $k_{i}^{(\pm r)}$, $e_{ij}^{(\pm r)}$, $f_{ji}^{(\pm r)}$ and $C$, and that the images of these monomials form a basis for $\mathrm{DY}_{h}(\mathfrak{gl}^{\mathfrak{s}}_{m|n})$.

For $~1\leq i<j\leq N, r\in\mathbb{Z}_{+}$, in $\widehat{\mathrm{DY}}_{h}(\mathfrak{gl}^{\mathfrak{s}}_{m|n})$, we define  inductively
\begin{align*}
e_{i,i+1}^{(\pm r)}=e_{i}^{(\pm r)}, \quad f_{i+1,i}^{(\pm r)}=f_{i}^{(\pm r)},\quad e_{ij}^{(\pm r)}=d_{j-1}\left[e_{i,j-1}^{(\pm r)},e_{j-1}^{(1)}\right],\quad~
f_{ji}^{(\pm r)}=d_{j-1}\left[f_{j-1}^{(1)},f_{j-1,i}^{(\pm r)}\right].
\end{align*}
These relations are clearly consistent with those in $\mathrm{DY}_{h}(\mathfrak{gl}^{\mathfrak{s}}_{m|n})$.

Let $\mathcal{E}$, $\mathcal{F}$ and $\mathcal{H}$ denote the subsuperalgebras of $\widehat{\mathrm{DY}}_{h}(\mathfrak{gl}^{\mathfrak{s}}_{m|n})$ generated by the elements $e_{j}^{(\pm r)}$, $f_{j}^{(\pm r)}$ and $k_{i}^{(\pm r)}$, respectively. Similarly, let $\mathcal{E}^{+}$, $\mathcal{F}^{+}$ and $\mathcal{H}^{+}$ be the sub-superalgebras generated by the positive elements $e_{j}^{(r)}$, $f_{j}^{(r)}$ and $k_{i}^{(r)}$, respectively, and let $\mathcal{E}^{-}$, $\mathcal{F}^{-}$ and $\mathcal{H}^{-}$ be those generated by the negative elements $e_{j}^{(-r)}$, $f_{j}^{(-r)}$ and $k_{i}^{(-r)}$, respectively. Define an ascending filtration on $\mathcal{E}$ by setting $\deg e_{i}^{(r)}=r-1,~ \deg
e_{i}^{(-r)}=-r$ and $\deg h=0$, and let $\operatorname{gr} \mathcal{E}$ be the associated graded superalgebra. Denote by $\bar{e}_{ij}^{(\pm r)}$ the image of $e_{ij}^{(\pm r)}$ in the corresponding graded component. As noted in \cite{G07}, $\mathcal{E}^{+}$ is spanned by ordered monomials in $e_{ij}^{(r)}$. The analogous spanning property for $\mathcal{E}$ follows from the relations
\begin{align*}
\left[\bar{e}_{ij}^{(r)},\bar{e}_{kl}^{(s)}\right]=d_{j}\delta_{kj}
\bar{e}_{il}^{(r+s-1)}-\theta_{i,j,k}\delta_{il}\bar{e}_{kj}^{(r+s-1)},\quad~\left[\bar{e}_{ij}^{(-r)},\bar{e}_{kl}^{(-s)}\right]=d_{j}\delta_{kj}
\bar{e}_{il}^{(-r-s)}-\theta_{i,j,k}\delta_{il}\bar{e}_{kj}^{(-r-s)},
\end{align*}
and
\begin{align*}
\left[\bar{e}_{ij}^{(r)},\bar{e}_{kl}^{(-s)}\right]&=
\left\{  
\begin{array}{ll}   
\delta_{kj}d_{j}
\bar{e}_{il}^{(r-s-1)}-\delta_{il}\theta_{i,j,k}\bar{e}_{kj}
^{(r-s-1)}, & r\leq s; \\
~\\    
\delta_{kj}d_{j}
\bar{e}_{il}^{(r-s)}-\delta_{il}\theta_{i,j,k}\bar{e}_{kj}
^{(r-s)}, & r>s.  
\end{array}
\right.
\end{align*}
Using \eqref{X1X1+}, \eqref{X2X1+}, and \eqref{e1eN-1}, these commutation relations can be verified in the same way as in \cite[Theorem 3]{G07}. Consequently, the algebra $\mathcal{E}$ is spanned by the ordered monomials in the elements $e_{ij}^{(r)}$ and $e_{kl}^{(-s)}$, taken in some fixed order. The same reasoning applies to $\mathcal{F}$.  For $\mathcal{H}$, the defining relations \eqref{kikj}–-\eqref{ki+ki-} imply that it is also spanned by ordered monomials in the $k_{i}^{(\pm r)}$. Furthermore, the defining relations of $\widehat{\mathrm{DY}_{h}}(\mathfrak{gl}^{\mathfrak{s}}_{m|n})$ implies that the multiplication map
\begin{align*}
\mathcal{F}~\widetilde{\otimes}~\mathcal{H}~\widetilde{\otimes}~
\mathcal{E}~\widetilde{\otimes}~\mathcal{A}[C]\rightarrow \widehat{\mathrm{DY}}_{h}(\mathfrak{gl}^{\mathfrak{s}}_{m|n})
\end{align*}
is surjective, where $\widetilde{\otimes}$ denotes the $h$-adic completion of the tensor product over $\mathbb{C}[[h]]$. Hence $\widehat{\mathrm{DY}}_{h}(\mathfrak{gl}^{\mathfrak{s}}_{m|n})$ is spanned by ordered monomials of the form
$$
\prod f_{ji}^{(\pm r)} \prod k_{i}^{(\pm r)}\prod e_{ij}^{(\pm r)} C^{a}, a\in\mathbb{Z}_{+}.
$$

On the other hand, in $\mathrm{DY}_{h}(\mathfrak{gl}^{\mathfrak{s}}_{m|n})$, applying Proposition \ref{Iso 1} to the matrix $T^{+}(u)$, we see that under the isomorphism \eqref{map1}, the images of $k_{i}^{(r)}$, $e_{ij}^{(r)}$ and $f_{ji}^{(r)}$ in the $(r-1)$-th graded component of gr$\mathrm{DY}_{h}(\mathfrak{gl}^{\mathfrak{s}}_{m|n})$  correspond to $-d_{i}E_{ii}^{(r-1)}$, $-d_{i}E_{ij}^{(r-1)}$ and $-d_{j}E_{ji}^{(r-1)}$, respectively. Similarly, the images of $k_{i}^{(-r)}$, $e_{ij}^{(-r)}$ and $f_{ji}^{(-r)}$ in the $(-r)$-th graded component correspond to $d_{i}E_{ii}^{(-r)}$, $d_{i}E_{ij}^{(-r)}$ and $d_{j}E_{ji}^{(-r)}$, respectively. By the PBW theorem for $U(\widehat{\mathfrak{gl}}^{\mathfrak{s}}_{m|n})$, the same ordered monomials form a basis of $\mathrm{DY}_{h}(\mathfrak{gl}^{\mathfrak{s}}_{m|n})$. It follows that $\phi$ is injective.
\end{proof}

We call $\widehat{\mathrm{DY}}_{h}(\mathfrak{gl}_{m|n}^{\mathfrak{s}})$ the Drinfeld presentation for the double super Yangian of $\mathfrak{gl}_{m|n}^{\mathfrak{s}}$. The same convention as in Remark \ref{conven-R} applies to $\widehat{\mathrm{DY}}_{h}(\mathfrak{gl}_{m|n}^{\mathfrak{s}})$.

\medskip
\section{Quantum Berezinian for the double super Yangian}

In this section, when necessary, we will add an additional superscript $\mathfrak{s}$ to our notation to avoid any ambiguity as $\mathfrak{s}$ varies.
The quantum Berezinian for the double super Yangian  $\mathrm{DY}^{0}_{-h}(\mathfrak{gl}_{m|n}^{\mathfrak{s}^{st}})$ was studied in \cite{BK25-1}. We now extend the construction to an arbitrary parity sequence $\mathfrak{s}$ by adapting the approach of \cite{CH23} to the super Yangian. We begin with the special case $m=n=1$.

\begin{lemm}
The following relations hold in $\mathrm{DY}^{0}_{h}(\mathfrak{gl}_{1|1}^{01})[[u^{\pm 1}]]:$
\begin{align}\label{e1}
t_{11}^{\pm}(u)t_{21}^{\pm}(u+h)&=t_{21}^{\pm}(u)t_{11}^{\pm}(u+h),\\ \label{e2}
t_{22}^{\pm}(u)t_{21}^{\pm}(u+h)&=t_{21}^{\pm}(u)t_{22}^{\pm}(u+h),\\ \label{e3}
t_{21}^{\pm}(u)t_{21}^{\pm}(u+h)&=0,\\ \label{e4}
t_{11}^{\pm}(u)t_{22}^{\pm}(u+h)-t_{22}^{\pm}(u)t_{11}^{\pm}(u+h)&=t_{12}^{\pm}(u)
t_{21}^{\pm}(u+h)+t_{21}^{\pm}(u)t_{12}^{\pm}(u+h).
\end{align}
\end{lemm}
\begin{proof}
The defining relations \eqref{Y1} and \eqref{Y2} yield the following identities:
\begin{align}\label{e5}
\left[t_{11}^{\pm}(u),t_{21}^{\pm}(u+h)\right]&=t_{21}^{\pm}(u)t_{11}^{\pm}(u+h)-t_{21}^{\pm}(u+h)t_{11}^{\pm}(u),\\ \label{e6}
\left[t_{22}^{\pm}(u),t_{21}^{\pm}(u+h)\right]&=t_{22}^{\pm}(u+h)t_{21}^{\pm}(u)-t_{22}^{\pm}(u)t_{21}^{\pm}(u+h),\\
\label{e7}
\left[t_{21}^{\pm}(u),t_{22}^{\pm}(u+h)\right]&=t_{21}^{\pm}(u+h)t_{22}^{\pm}(u)-t_{21}^{\pm}(u)t_{22}^{\pm}(u+h),\\
\label{e8}
\left[t_{21}^{\pm}(u),t_{21}^{\pm}(u+h)\right]&=t_{21}^{\pm}(u+h)t_{21}^{\pm}(u)-t_{21}^{\pm}(u)t_{21}^{\pm}(u+h),\\
\label{e9}
\left[t_{11}^{\pm}(u),t_{22}^{\pm}(u+h)\right]&=t_{21}^{\pm}(u)t_{12}^{\pm}(u+h)-t_{21}^{\pm}(u+h)t_{12}^{\pm}(u),\\
\label{e10}
\left[t_{12}^{\pm}(u),t_{21}^{\pm}(u+h)\right]&=t_{12}^{\pm}(u)t_{21}^{\pm}(u+h)+t_{21}^{\pm}(u+h)t_{12}^{\pm}(u).
\end{align}
By \eqref{e5}, we can deduce \eqref{e1}. Using \eqref{e6} and \eqref{e7}, we obtain
\begin{align*}
2~t_{22}^{\pm}(u)t_{21}^{\pm}(u+h)=t_{22}^{\pm}(u+h)t_{21}^{\pm}(u)+t_{21}^{\pm}(u+h)
t_{22}^{\pm}(u)=2~t_{21}^{\pm}(u)t_{22}^{\pm}(u+h).
\end{align*}
This proves \eqref{e2}.
From \eqref{e8}, we get
\begin{align*}
t_{21}^{\pm}(u)t_{21}^{\pm}(u+h)+t_{21}^{\pm}(u+h)t_{21}^{\pm}(u)=t_{21}^{\pm}(u+h)t_{21}^{\pm}(u)-t_{21}^{\pm}(u)t_{21}^{\pm}(u+h),
\end{align*}
so $2~t_{21}^{\pm}(u)t_{21}^{\pm}(u+h)=0$, and \eqref{e3} follows directly. According to \eqref{e9} and \eqref{e10}, we have
\begin{align*}
t_{11}^{\pm}(u)t_{22}^{\pm}(u+h)-t_{22}^{\pm}(u)t_{11}^{\pm}(u+h)=t_{12}^{\pm}(u)t_{21}^{\pm}(u+h)+t_{21}^{\pm}(u)t_{12}^{\pm}(u+h).
\end{align*}
Then relation \eqref{e4} holds.
\end{proof}

By Lemma \ref{DY-iso}, the definition of $\mathrm{DY}^{0}_{h}(\mathfrak{gl}_{m|n}^{\mathfrak{s}})$ is independent of the choice of $\mathfrak{s}$. In particular, the double Yangian $\mathrm{DY}^{0}_{h}(\mathfrak{gl}_{1|1}^{10})$ and $\mathrm{DY}^{0}_{h}(\mathfrak{gl}_{1|1}^{01})$ are isomorphic via the map $t_{ij}^{\pm}(u)\mapsto t_{3-i,3-j}^{\pm}(u)$.

\begin{lemm} \label{map DY10}
Let $\sigma_{1}: \mathrm{DY}^{0}_{h}(\mathfrak{gl}_{1|1}^{10})\rightarrow \mathrm{DY}^{0}_{h}(\mathfrak{gl}_{1|1}^{01})$ be the superalgebra isomorphism defined previously. Then $\sigma_1$ sends
\begin{align*}
k_{1}^{\pm }(u)^{-1}k_{2}^{\pm }(u)\mapsto k_{1}^{\pm }(u+h)k_{2}^{\pm }(u+h)^{-1}.
\end{align*}
\end{lemm}
\begin{proof}
By the definition of quasideterminants, we have in $\mathrm{DY}^{0}_{h}(\mathfrak{gl}_{1|1}^{10})$,
\begin{align*}
k_{1}^{\pm }(u)^{-1}k_{2}^{\pm }(u)=t_{11}^{\pm }(u)^{-1}\Big(t_{22}^{\pm }(u)-t_{21}^{\pm }(u)t_{11}^{\pm }(u)^{-1}t_{12}^{\pm }(u)\Big);
\end{align*}
hence, $\sigma_1$ sends
\begin{align*}
k_{1}^{\pm}(u)^{-1}k_{2}^{\pm}(u)\mapsto t_{22}^{\pm }(u)^{-1}
\Big(t_{11}^{\pm }(u)-t_{12}^{\pm }(u)t_{22}^{\pm}(u)^{-1}t_{21}^{\pm }(u)\Big).
\end{align*}

It remains to verify in $\mathrm{DY}^{0}_{h}(\mathfrak{gl}_{1|1}^{01})$,
\begin{gather*}
    t_{22}^{\pm}(u)t_{11}^{\pm }(u+h)=\Big(t_{11}^{\pm}(u)-t_{12}^{\pm}(u)t_{22}^{\pm}(u)^{-1}t_{21}^{\pm}(u)\Big)k_{2}^{\pm}(u+h),
\end{gather*}
which follows directly from \eqref{e1}-\eqref{e4}.
\end{proof}

For a fixed $\mathfrak{s}\in\mathcal{S}_{m|n}$, we adopt the following notations from \cite{P16}:
\begin{itemize}
  \item $\check{\mathfrak{s}}$ := the $0^{n}1^{m}$-sequence obtained by interchanging the 0s and 1s of $\mathfrak{s}$.
  \item $\mathfrak{s}^{r}$ := the reverse of $\mathfrak{s}$.
  \item $\mathfrak{s}^{\dag}$ := $(\check{\mathfrak{s}})^{r}$, the reverse of $\check{\mathfrak{s}}$.
\end{itemize}
For example, if $\mathfrak{s}=101101100$, then $\check{\mathfrak{s}}=010010011$, $\mathfrak{s}^{r}=001101101$, and $\mathfrak{s}^{\dag}=110010010$. Moreover, the concatenation of two parity sequences $\mathfrak{s}$ and $\mathfrak{s}'$ is simply denoted by $\mathfrak{s}\mathfrak{s}'$.

The following results generalize the corresponding results for the super Yangian (cf. \cite{CH23, G07, P16}).

\begin{lemm}\label{Hom 1}
There exist a superalgebra isomorphisms
$$
\omega_{m|n}:\mathrm{DY}^{0}_{h}(\mathfrak{gl}_{m|n}^{\mathfrak{s}})\rightarrow\mathrm{DY}^{0}_{h}
(\mathfrak{gl}_{m|n}^{\mathfrak{s}})
$$
such that for all $i,j\in I_{\mathfrak{s}}$,
$$
\omega_{m|n}(T^{\pm}(u))=T^{\pm}(-u)^{-1}.
$$
\end{lemm}
\begin{proof}
One checks directly that $\rho_{m|n}$ and $\omega_{m|n}$ preserve the relation \eqref{RTT1}.
\end{proof}

Let $p,q\in\mathbb{N}$ and let $\mathfrak{s}'$ be an arbitrary $0^{p}1^{q}$-sequence. Denote by
$$
\varphi_{p|q}:\mathrm{DY}^{0}_{h}(\mathfrak{gl}_{m|n}^{\mathfrak{s}})\hookrightarrow
\mathrm{DY}^{0}_{h}(\mathfrak{gl}_{m+p|n+q}^{\mathfrak{s}'\mathfrak{s}}).
$$
The injective homomorphism $\varphi_{p|q}$ is defined by the assignment $t_{ij}^{(\pm r)}\mapsto t_{i+p+q,j+p+q}^{(\pm r)}$. Furthermore, define $\psi_{p|q}:\mathrm{DY}^{0}_{h}(\mathfrak{gl}_{m|n}^{\mathfrak{s}})\hookrightarrow
\mathrm{DY}^{0}_{h}(\mathfrak{gl}_{m+p|n+q}^{\mathfrak{s}'\mathfrak{s}})$ to be the injective homomorphism defined as the composition
$$
\psi_{p|q}=\omega_{p+m|q+n}\circ\varphi_{p|q}\circ\omega_{m|n}.
$$

In view of \cite[Lemma 4.1]{P16}, we have the following lemma.
\begin{lemm}\label{Hom 3}
For any $p,q\in\mathbb{N}$, we have
\begin{equation} \label{psi}
\psi_{p|q}(t^{\pm}_{ij}(u))=\left|\begin{array}{cccc}
	t^{\pm}_{11}(u)&\cdots&t^{\pm}_{1,p+q}(u)&t^{\pm}_{1,p+q+j}(u)\\
	\vdots&&\vdots&\vdots\\
	t^{\pm}_{p+q,1}(u)&\cdots&t^{\pm}_{p+q,p+q}(u)&t^{\pm}_{p+q,p+q+j}(u)\\ t^{\pm}_{p+q+i,1}(u)&\cdots&t^{\pm}_{p+q+i,p+q}(u)&\framebox{$t^{\pm}_{p+q+i,p+q+j}(u)$}\\
\end{array}\right|.
\end{equation}
\end{lemm}

As an immediate consequence of Lemma \ref{Hom 3}, we have the following lemma.

\begin{lemm}\label{Hom 4}
For all $i\in I_{\mathfrak{s}}$ and $j\in I_{\mathfrak{s}}'$, the mapping $\psi_{p|q}$ acts on the Gaussian generators by
$$
k_{i}^{\pm}(u)\mapsto k_{p+q+i}^{\pm}(u),\quad e_{j}^{\pm}(u)\mapsto e_{p+q+j}^{\pm}(u),\quad f_{j}^{\pm}(u)\mapsto f_{p+q+j}^{\pm}(u).
$$
\end{lemm}

\begin{rema}
By \eqref{psi}, $\psi_{p|q}$ depends only on $p+q$, so we write $\psi_{p|q}=\psi_{p+q}$ when appropriate.
\end{rema}

For $\mathfrak{s}=(s_1 s_{2} \cdots s_{N})\neq \mathfrak{s}^{st}$, choose an index $i$ with $(s_{i}s_{i+1})=(10)$, and let $\bar{\mathfrak{s}}$ be the sequence obtained from $\mathfrak{s}$ by switching $s_{i}$ and $s_{i+1}$. The simple reflection $\sigma_{i}=(i,i+1)$ for $i\in I_{\mathfrak{s}}'$ induces an isomorphism:
\begin{align*}
\sigma_{i}: \mathrm{DY}^{0}_{h}(\mathfrak{gl}_{m|n}^{\mathfrak{s}})\rightarrow \mathrm{DY}^{0}_{h}(\mathfrak{gl}_{m|n}^{\bar{\mathfrak{s}}});\quad~
 t_{kl}^{\pm }(u)\mapsto t_{\sigma_{i}(k),\sigma_{i}(l)}^{\pm}(u).
\end{align*}

\begin{prop} \label{sigma}
The mapping $\sigma_i$ sends $k_{i}^{\pm }(u)^{-1}k_{i+1}^{\pm }(u)\mapsto k_{i}^{\pm}(u+h)k_{i+1}^{\pm}(u+h)^{-1}$.
\end{prop}
\begin{proof}
If $i=1$, then $(s_1s_2)=(10)$. There is a standard embedding $\mathrm{DY}^{0}_{h}
(\mathfrak{gl}_{1|1}^{10})\hookrightarrow\mathrm{DY}^{0}_{h}(\mathfrak{gl}_{m|n}^{\mathfrak{s}})$, given by $t_{ij}^{\pm}(u)\mapsto t_{ij}^{\pm}(u)$. The claim then follows immediately from Lemma \ref{map DY10}.

Suppose now that $i>1$. We decompose $\mathfrak{s}=\mathfrak{s}_{1} \mathfrak{s}_{2}$ with $\mathfrak{s}_{1}:=(s_1\cdots s_{i-1})$ and $\mathfrak{s}_{2}:=(s_{i}\cdots s_{N})$. Let $\bar{\mathfrak{s}}_{2}$
  denote the sequence obtained from $\mathfrak{s}_{2}$ by swapping $s_{i}$ and $s_{i+1}$, so $\bar{\mathfrak{s}}=\mathfrak{s}_{1}\bar{\mathfrak{s}}_{2}$. Applying the shift map \eqref{psi}, together with Lemma \ref{Hom 4}, we obtain in $\mathrm{DY}^{0}_{h}(\mathfrak{gl}_{m|n}^{\mathfrak{s}})$,
\begin{align*}
\sigma_{i}\left(k_{i}^{\pm}(u)^{-1}k_{i+1}^{\pm }(u)\right)&=\sigma_{i}\left(\psi_{i-1}\left(k_{1}^{\pm }(u)^{-1}k_{2}^{\pm }(u)\right)\right)
=\psi_{i-1}\left(\sigma_{1}\left(k_{1}^{\pm }(u)^{-1}k_{2}^{\pm }(u)\right)\right).
\end{align*}
In $\mathrm{DY}^{0}_{h}(\mathfrak{gl}_{m|n}^{\bar{\mathfrak{s}}})$, we have
\begin{align*}
k_{i}^{\pm }(u+h)k_{i+1}^{\pm }(u+h)^{-1}=\psi_{i-1}\left(k_{1}^{\pm }(u+h)k_{2}^{\pm }(u+h)^{-1}\right).
\end{align*}
Note that $\sigma_{1}$ sends $k_{1}^{\pm }(u)^{-1}k_{2}^{\pm }(u)\mapsto k_{1}^{\pm}(u+h)k_{2}^{\pm}(u+h)^{-1}$, which completes the proof.
\end{proof}

Given a 01-sequence $\mathfrak{s}=(s_{1}\cdots s_{N})$, there exist indices $1\leq i_{1}< i_{2}<\cdots<i_{m}\leq N$ and $1\leq j_{1}< j_{2}<\cdots<j_{n}\leq N$ such that $s_{i_{k}}=0$ and $s_{j_{l}}=1$ for all $1\leq k\leq m, 1\leq l\leq n$. Obviously, the indices $i_k$ and $j_l$ are uniquely determined by $\mathfrak{s}$.

\begin{defi}\label{quantum berezinian}
The quantum Berezinian for the matrix $T^{\pm}(u)$ of $\mathrm{DY}^{0}_{h}(\mathfrak{gl}_{m|n}^{\mathfrak{s}})$ is defined via
\begin{equation}\label{berezinian 1}
    \begin{split}
        b_{m|n}^{\pm (\mathfrak{s})}(u)=&\sum_{\rho}\operatorname{sgn}(\rho)t_{i_{\rho(1)},i_{1}}^{\pm }(u)t_{i_{\rho(2)},i_{2}}^{\pm }(u+h)\cdots t_{i_{\rho(m)},i_{m}}^{\pm }(u+(m-1)h)\\ 
        &~~\times \sum_{\sigma}\operatorname{sgn}(\sigma){\tilde{t}_{j_{1},j_{\sigma(1)}}^{\pm }}(u+(m-1)h)\cdots {\tilde{t}_{j_{n},j_{\sigma(n)}}^{\pm }}(u+(m-n)h),
    \end{split}
\end{equation}
where $\rho$ and $\sigma$ range over all permutations of orders $m$ and $n$, respectively, and $T^{\pm}(u)^{-1}:=
\left(\tilde{t}_{ij}^{\pm}(u)\right)_{i,j=1}^{N}$.
\end{defi}

\begin{rema}
When $\mathfrak{s}=\mathfrak{s}^{st}$ is the standard $0^{m}1^{n}$-sequence, we have $i_{k}=k$ and $j_{l}=m+l$. In this case, \eqref{berezinian 1} is consistent with the definition in \cite{BK25-1} after replacing $h$ by $-h$.
\end{rema}

Observe that the permutation $\sigma_{\mathfrak{s}}$ of order $N$ sending $i_{k}\mapsto k, j_{l}\mapsto m+l$, induces an isomorphism (still denoted $\sigma_{\mathfrak{s}}$).
\begin{align*}
\sigma_{\mathfrak{s}}: \mathrm{DY}^{0}_{h}(\mathfrak{gl}_{m|n}^{\mathfrak{s}})\rightarrow \mathrm{DY}^{0}_{h}(\mathfrak{gl}_{m|n}^{\mathfrak{s}^{st}});\quad~
t_{ij}^{\pm }(u)\mapsto t_{\sigma_{\mathfrak{s}}(i),\sigma_{\mathfrak{s}}(j)}^{\pm }(u).
\end{align*}
In particular, we have $\sigma_{\mathfrak{s}}\left(b_{m|n}^{\pm (\mathfrak{s})}(u)\right)=b_{m|n}^{\pm (\mathfrak{s}^{st})}(u)$. 

\begin{lemm}\label{berezinian 2}
The quantum Berezinian can be expressed in terms of the Drinfeld generators as follows
\begin{align}\label{berezinian Drinfeld}
b_{m|n}^{\pm (\mathfrak{s})}(u)=\prod_{i=1}^{N}\tilde{k}_{i}^{\pm}(u_i),\quad
 {\rm where~}
 \tilde{k}_{i}^{\pm }(u):=
\left\{\begin{array}{ll}
k_{i}^{\pm}(u),  & {\rm if~} s_{i}=0, \\
k_{i}^{\pm }(u)^{-1}, & {\rm if~} s_{i}=1,
\end{array}
\right.
\end{align}
while $u_{1}=u-s_{1}h$ and $u_{i+1}=\left\{\begin{array}{ll}
u_{i}+h, & {\rm if~} s_{i}=s_{i+1}=0, \\
 u_{i}-h, & {\rm if~} s_{i}=s_{i+1}=1, \\
 u_{i}, & {\rm if~} s_{i}\neq s_{i+1}.
\end{array}
 \right.$
\end{lemm}
\begin{proof}
If $\mathfrak{s}=\mathfrak{s}^{st}$, the assertion follows directly from \cite[Lemma 3.6]{BK25-1}. Otherwise, let $l\in I_{\mathfrak{s}}'$ be the smallest index such that $(s_{l}s_{l+1})=(10)$. Denote by $\sigma_{l}$ the isomorphism induced by the simple reflection $(l,l+1)$, and let $\bar{\mathfrak{s}}:=(\bar{s}_{1}\cdots\bar{s}_{N})$ be the sequence obtained from $\mathfrak{s}$ by switching $s_{l}$ and $s_{l+1}$. In \eqref{berezinian Drinfeld}, we have $\tilde{k}_{l}(u_{l})\tilde{k}_{l+1}(u_{l+1})=k_{l}(u_{l})^{-1}k_{l+1}(u_{l})$. Proposition \ref{sigma} then gives $\sigma_{l}(b_{m|n}^{\pm (\mathfrak{s})}(u))=b_{m|n}^{\pm (\bar{\mathfrak{s}})}(u)$. Repeating this procedure until the standard sequence is reached, we obtain a finite sequence of simple reflections whose composition is the permutation $\sigma_{\mathfrak{s}}$, hence $\sigma_{\mathfrak{s}}(b_{m|n}^{\pm (\mathfrak{s})}(u))=b_{m|n}^{\pm (\mathfrak{s}^{st})}(u)$. The claim then follows from \cite[Lemma 3.6]{BK25-1}.
\end{proof}

Combining this with \cite[Proposition 3.7 and Theorem 3.10]{BK25-1}, we obtain the following results.

\begin{prop}
The coefficients of the quantum Berezinian lie in the center of~ $\mathrm{DY}^{0}_{h}(\mathfrak{gl}_{m|n}^{\mathfrak{s}})$. 
\end{prop}

Let $b^{(\pm r)}_{m|n}$ for $r\in\mathbb{Z}_{+}$ denote the coefficients of the quantum Berezinian, so that
\begin{align*}
b^{+}_{m|n}(u)=1-h\sum_{r\geqslant 1}b^{(r)}_{m|n}u^{-r}\quad{\rm and}~\quad b^{-}_{m|n}(u)=1+h\sum_{r\geqslant 1}b^{(-r)}_{m|n}u^{r-1}.
\end{align*}

\begin{coro}
The coefficients $b^{(\pm r)}_{m|n}$ with $r\in\mathbb{Z}_{+}$ form a set of algebraically independent topological generators of the center of $\mathrm{DY}^{0}_{h}(\mathfrak{gl}_{m|n}^{\mathfrak{s}})$.
\end{coro}

\medskip
\section{Double super Yangian of the special linear Lie superalgebra}

In this section, for any $\mathfrak{s}\in\mathcal{S}_{m|n}$, we first use the quantum Berezinian of $\mathrm{DY}^{0}_{h}(\mathfrak{gl}_{m|n}^{\mathfrak{s}})$ (defined in the previous section) to introduce the R-matrix presentation of the double Yangian for the Lie superalgebra $\mathfrak{sl}_{m|n}^{\mathfrak{s}}$. We then give its Drinfeld presentation and finally prove that these two presentations are isomorphic.

\begin{defi}
The double Yangian $\mathrm{DY}_{h}(\mathfrak{sl}^{\mathfrak{s}}_{m|n})$ associated with the special linear Lie superalgebra $\mathfrak{sl}^{\mathfrak{s}}_{m|n}$ is defined as the following subalgebra of $\mathrm{DY}_{h}(\mathfrak{gl}^{\mathfrak{s}}_{m|n}):$
$$
\mathrm{DY}_{h}(\mathfrak{sl}^{\mathfrak{s}}_{m|n}):=\{y\in \mathrm{DY}_{h}(\mathfrak{gl}^{\mathfrak{s}}_{m|n})\mid \mu_{g^{\pm}}(y)=y \text{ for all } g^{\pm}\}.
$$
Here the map $\mu_{g^{\pm}}$ is the automorphism of $\mathrm{DY}_{h}(\mathfrak{gl}^{\mathfrak{s}}_{m|n})$ given by
$$
\mu_{g^{\pm}}: T^{\pm}(u)\mapsto g^{\pm}(u)T^{\pm}(u),\quad~C\mapsto C,
$$
where the formal series $g^{+}(u)=1+h\sum_{r\geqslant 1}g_{r}^{+}u^{-r}$ and $g^{-}(u)=1+h\sum_{r\geqslant 1}g_{r}^{-}u^{r-1}$ lie in $\mathcal{A}[[u^{\mp 1}]]$.
\end{defi}

\begin{lemm}\label{mug_action}
For any $1\leqslant i, j\leqslant N$ with $i<j$, we have
\begin{align}\label{ki}
\mu_{g^{\pm}}(k_{i}^{\pm}(u))&=g^{\pm}(u)k_{i}^{\pm}(u);\\
\mu_{g^{\pm}}(e_{ij}^{\pm}(u))&=e_{ij}^{\pm}(u);\\ \label{fij}
\mu_{g^{\pm}}(f_{ji}^{\pm}(u))&=f_{ji}^{\pm}(u).
\end{align}
\end{lemm}
\begin{proof}
The proof is similar to that of \cite[Lemma 4.4]{XLZ}.
\end{proof}

Let $\mathrm{DY}_{h}^{0}(\mathfrak{sl}_{m|n}^{\mathfrak{s}}):=\mathrm{DY}_{h}(\mathfrak{sl}_{m|n}^{\mathfrak{s}})/\langle C\rangle$ be the level-0 double Yangian of $\mathfrak{sl}_{m|n}^{\mathfrak{s}}$, obtained by setting the central element $C$ to zero. The superalgebra $\mathrm{DY}_{h}(\mathfrak{sl}_{m|n}^{\mathfrak{s}})$ is a central extension of $\mathrm{DY}_{h}^{0}(\mathfrak{sl}_{m|n}^{\mathfrak{s}})$ by adding the central generator $C$, giving the isomorphism
\begin{align*}
\mathrm{DY}_{h}(\mathfrak{sl}_{m|n}^{\mathfrak{s}})\cong \mathrm{DY}_{h}^{0}(\mathfrak{sl}_{m|n}^{\mathfrak{s}})\otimes_{\mathcal{A}}\mathcal{A}[C].
\end{align*}

\begin{prop}
Let $Z^{\mathfrak{s}}_{m|n}$ be the subalgebra of $\mathrm{DY}^{0}_{h}(\mathfrak{gl}^{\mathfrak{s}}_{m|n})$ generated by the coefficients of the quantum Berezinian. Then for $m\neq n$, we have
$$
\mathrm{DY}_{h}(\mathfrak{gl}^{\mathfrak{s}}_{m|n})\cong Z^{\mathfrak{s}}_{m|n} \otimes \mathrm{DY}_{h}(\mathfrak{sl}^{\mathfrak{s}}_{m|n}).
$$
\end{prop}
\begin{proof}
We assume that $m>n$. (The result for $n>m$ follows from this by applying the map $\zeta_{m|n}$.) By Proposition 2.15  of \cite{MNO96}, there exists a unique series $\bar{b}_{m|n}^{\pm }(u)\in\mathrm{DY}^{0}_{h}(\mathfrak{gl}^{\mathfrak{s}}_{m|n})[[u^{\mp}]]$ such that
$$
b_{m|n}^{\pm (\mathfrak{s})}(u)=\bar{b}_{m|n}^{\pm }(u)\bar{b}_{m|n}^{\pm}(u+h)\cdots \bar{b}_{m|n}^{\pm}(u+(m-n-1)h).
$$
From Lemma \ref{berezinian 2} and \eqref{ki}, we have
$$
\mu_{g^{\pm}}\left(b_{m|n}^{\pm (\mathfrak{s})}(u)\right)=g^{\pm}(u)g^{\pm}(u+h)\cdots g^{\pm}(u+(m-n-1)h)b_{m|n}^{\pm (\mathfrak{s})}(u).
$$
The uniqueness of this factorization forces $\mu_{g^{\pm}}\left(\bar{b}_{m|n}^{\pm}(u)\right)=g^{\pm}(u) \bar{b}_{m|n}^{\pm}(u)$. Now, we define $\bar{t}_{ij}^{\pm}(u)=\left(\bar{b}_{m|n}^{\pm}(u)
\right)^{-1}t_{ij}^{\pm}(u)$. Then $\mu_{g^{\pm}}\left(\bar{t}_{ij}^{\pm}(u)
\right)=\bar{t}_{ij}^{\pm}(u)$, so $\bar{t}_{ij}^{\pm}(u)\in \mathrm{DY}^{0}_{h}(\mathfrak{sl}^{\mathfrak{s}}_{m|n})$. This implies that
$$
\mathrm{DY}^{0}_{h}(\mathfrak{gl}^{\mathfrak{s}}_{m|n})= Z^{\mathfrak{s}}_{m|n} \cdot \mathrm{DY}^{0}_{h}(\mathfrak{sl}^{\mathfrak{s}}_{m|n}).
$$
Moreover,
$Z^{\mathfrak{s}}_{m|n} \cap\mathrm{DY}^{0}_{h}(\mathfrak{sl}^{\mathfrak{s}}_{m|n})$ is trivial, hence $\mathrm{DY}^{0}_{h}(\mathfrak{gl}^{\mathfrak{s}}_{m|n})\cong Z^{\mathfrak{s}}_{m|n} \otimes \mathrm{DY}^{0}_{h}(\mathfrak{sl}^{\mathfrak{s}}_{m|n})$. Consequently, 
\begin{align*}
\mathrm{DY}_{h}(\mathfrak{gl}^{\mathfrak{s}}_{m|n})\cong\mathrm{DY}^{0}_{h}(\mathfrak{gl}^{\mathfrak{s}}_{m|n})\otimes_{\mathcal{A}}\mathcal{A}[C]\cong Z^{\mathfrak{s}}_{m|n} \otimes \mathrm{DY}^{0}_{h}(\mathfrak{sl}^{\mathfrak{s}}_{m|n})\otimes_{\mathcal{A}}\mathcal{A}[C]\cong Z^{\mathfrak{s}}_{m|n} \otimes \mathrm{DY}_{h}(\mathfrak{sl}^{\mathfrak{s}}_{m|n}).
\end{align*}
\end{proof}

The following lemma gives the R-matrix presentation of $\mathrm{DY}_{h}(\mathfrak{sl}^{\mathfrak{s}}_{m|n})$. Its proof is analogous to those of \cite[Lemma 4.7]{XLZ} and \cite[Lemma 7]{G07}.
\begin{lemm} \label{SL RTT}
For any $m,n\geqslant 0$, the coefficients of the series
\begin{equation}\label{SL}
{k_{i}^{\pm}(u)}^{-1}k_{i+1}^{\pm}(u),\quad e_{i}^{\pm}(u),\quad f_{i}^{\pm}(u), \quad {\rm for~} i\in I_{\mathfrak{s}}',
\end{equation}
along with the central element $C$ generate the subalgebra $\mathrm{DY}_{h}(\mathfrak{sl}^{\mathfrak{s}}_{m|n})$.
\end{lemm}
\begin{proof}
The algebra $\mathrm{DY}_{h}(\mathfrak{gl}^{\mathfrak{s}}_{m|n})$ is generated by the coefficients of the series $\{k_{j}^{\pm}(u),e_i^{\pm}(u), f_{i}^{\pm}(u)\mid j\in I_{\mathfrak{s}}, i\in I_{\mathfrak{s}}'\}$ together with the central element $C$, so the coefficients of the series $k_{1}^{\pm}(u)$ together with those in \eqref{SL} also generate the algebra $\mathrm{DY}_{h}(\mathfrak{gl}^{\mathfrak{s}}_{m|n})$. Note that, for any $g^{\pm}$, the automorphism $\mu_{g^{\pm}}$ fixes all generators in (\ref{SL}) and $\mu_{g^{\pm}}(k_{1}^{\pm}(u))=g^{\pm}(u)k_{1}^{\pm}(u)$.

By Theorem \ref{PBW basis} and \eqref{Gauss ii}-\eqref{Gauss ji}, any element $J\in \mathrm{DY}_{h}(\mathfrak{gl}^{\mathfrak{s}}_{m|n})$ is a polynomial in $k_{1}^{(\pm 1)}$, $k_{1}^{(\pm 2)}, \cdots$ and the other generators fixed by all $\mu_{g^{\pm}}$.  We may assume that the monomials in $J$ are ordered with $f_{j}^{(\pm r)}$ preceding $k_{i}^{(\pm r)}$, which precede $e_{j}^{(\pm r)}$. Suppose $J\in \mathrm{DY}_{h}(\mathfrak{sl}^{\mathfrak{s}}_{m|n})$. Let $M$ denote the maximal index $r$ for which the generator $k_{1}^{(\pm r)}$ appears in $J$, and let $H$ denote the highest power of any such $k_{1}^{(\pm r)}$ involved in $J$. We then write:
\begin{align*}
J=\sum_{a}F_a(k_{1}^{(1)})^{a_1^+}(k_{1}^{(-1)})^{a_1^-}\cdots
(k_{1}^{(M)})^{a_M^+}(k_{1}^{(-M)})^{a_M^-}K_aE_a,
\end{align*}
where $E_a$, $K_a$, $F_a$ are monomials in the generators fixed by $\mu_{g^{\pm}}$, and the sum is over all $2M$-tuples $a=(a_1^{+},a_1^{-},\cdots,a_M^{+},a_M^-)$ with $0\leqslant a_{i}^{\pm}\leqslant H$. Fix $g^{+}(u)=1+\lambda h u^{- M}$ and $g^{-}(u)=1+\lambda h u^{M-1}$ with $\lambda\in\mathbb{C}^{*}$. Then we have
\begin{align*}
\mu_{g^{\pm}}(J)=\sum_{a}F_a(k_{1}^{(1)})^{a_1^+}(k_{1}^{(-1)})^{a_1^-}\cdots(\lambda+k_{1}^{(M)})^{a_M^+}(\lambda+k_{1}^{(-M)})^{a_M^-}K_aE_a=J.
\end{align*}
Due to the linear independence of the different monomials and the arbitrariness of $\lambda$, it follows that $k_{1}^{(\pm M)}$ cannot appear in $J$. This implies the claim.
\end{proof}

For $m=n$, the coefficients of the quantum Berezinian lie in the center of $\mathrm{DY}_{h}(\mathfrak{sl}^{\mathfrak{s}}_{n|n})$. We define the double Yangian associated with the classical Lie superalgebra $A(n-1,n-1)$ as the quotient
$$
\mathrm{DY}_{h}(\mathfrak{psl}^{\mathfrak{s}}_{n|n}):=\mathrm{DY}_{h}(\mathfrak{sl}^{\mathfrak{s}}_{n|n}
)/\langle b_{n|n}^{\pm}(u)=1\rangle=\mathrm{DY}_{h}
(\mathfrak{sl}^{\mathfrak{s}}_{n|n})/Z^{\mathfrak{s}}_{n|n}.
$$

Next, we establish the Drinfeld presentation of $\mathrm{DY}_{h}(\mathfrak{sl}^{\mathfrak{s}}_{m|n})$. Introduce the currents:
\begin{align*}
H_{i}^{\pm}(u)&=k_{i+1}^{\pm}\left(u+\frac{1}{2}h\nu_{i}\right)k_{i}^{\pm}\left(u+\frac{1}{2}h\nu_{i}\right)^{-1},\\
E_{i}(u)&=\frac{(-1)^{|i||i+1|}}{h}X_{i}^{+}(u+\frac{1}{2}h\nu_{i}),\\
F_{i}(u)&=\frac{1}{h}X_{i}^{-}(u+\frac{1}{2}h\nu_{i}),
\end{align*}
where $\nu_{i}=\sum_{j=1}^{i}d_{j}$. By Lemma \ref{SL RTT}, the coefficients of the series $H_{i}^{\pm}(u)$, $E_{i}(u)$ and $F_{i}(u)$ for $i\in I_{\mathfrak{s}}'$, along with the central element $C$ generate the double super Yangian $\mathrm{DY}_{h}(\mathfrak{sl}_{m|n}^{\mathfrak{s}})$.

\begin{lemm}\label{Drinfeld generator sl}
The following relations hold in the double Yangian $\mathrm{DY}_{h}(\mathfrak{sl}^{\mathfrak{s}}_{m|n}):$
\begin{align}\label{HiHj}
&H_{i}^{\pm}(u)H_{j}^{\pm}(v)=H_{j}^{\pm}(v)H_{i}^{\pm}(u),\\ \label{Hi+Hj-}
&H_{i}^{\pm}(u)H_{j}^{\mp}(v)=\frac{(u_{\mp}-v_{\pm}-B_{ij}h)(u_{\pm}-v_{\mp}+B_{ij}h)}{(u_{\mp}-v_{\pm}+B_{ij}h)(u_{\pm}-v_{\mp}-B_{ij}h)}H_{j}^{\mp}(v)H_{i}^{\pm}(u),\\ \label{HiEj}
&H_{i}^{\pm}(u)^{-1}E_{j}(v)H_{i}^{\pm}(u)=\frac{u_{\pm}-v-B_{ij}h}{u_{\pm}-v+B_{ij}h}E_{j}(v),\\ \label{HiFj}
&H_{i}^{\pm}(u)F_{j}(v)H_{i}^{\pm}(u)^{-1}=\frac{u_{\mp}-v-B_{ij}h}{u_{\mp}-v+B_{ij}h}F_{j}(v),\\ \label{EiEj}
&E_{i}(u)E_{j}(v)=(-1)^{|\alpha_{i}||\alpha_{j}|}\frac{u-v+B_{ij}h}{u-v-B_{ij}h}E_{j}(v)E_{i}(u),\\ \label{FiFj}
&F_{i}(u)F_{j}(v)=(-1)^{|\alpha_{i}||\alpha_{j}|}\frac{u-v-B_{ij}h}{u-v+B_{ij}h}F_{j}(v)F_{i}(u),\\ \label{EiFj}
&\left[E_{i}(u),F_{j}(v)\right]=\frac{\delta_{ij}}{h}\{\delta(\frac{v_{+}}{u_{-}})H_{i}^{+}(u_{-})-\delta(\frac{v_{-}}{u_{+}})H_{i}^{-}(v_{-})\},\\ \label{serre E1}
&{\rm Sym}_{u_{1},u_{2}}\left[E_{i}(u_{1}),\left[E_{i}(u_{2}),E_{i\pm 1}(v)\right]\right]=0,\quad |\alpha_{i}|=0,  ~i\pm 1 \in I_{\mathfrak{s}}'\\ \label{serre F1}
&{\rm Sym}_{u_{1},u_{2}}\left[F_{i}(u_{1}),\left[F_{i}(u_{2}),F_{i\pm 1}(v)\right]\right]=0,\quad |\alpha_{i}|=0,  ~i\pm 1 \in I_{\mathfrak{s}}'\\ \label{serre E2}
&{\rm Sym}_{u_{1},u_{2}}\left[E_{i}(u_{1}),\left[E_{i+1}(v_{1}),\left[E_{i}(u_{2}),
E_{i-1}(v_{2})\right]\right]\right]=0,\quad |\alpha_{i}|=1, ~i\pm 1 \in I_{\mathfrak{s}}'\\ \label{serre F2}
&{\rm Sym}_{u_{1},u_{2}}\left[F_{i}(u_{1}),\left[F_{i+1}(v_{1}),\left[F_{i}(u_{2}),
F_{i-1}(v_{2})\right]\right]\right]=0,\quad |\alpha_{i}|=1  ~i\pm 1 \in I_{\mathfrak{s}}'.
\end{align}
Here, we set $B_{ij}=\frac{1}{2}a_{ij}$, where $A=(a_{ij})$ is the Cartan matrix of the Lie superalgebra $\mathfrak{sl}^{\mathfrak{s}}_{m|n}$.
\end{lemm}
\begin{proof}
From \eqref{kikj}, we obtain \eqref{HiHj}. Next, we prove that \eqref{Hi+Hj-}-\eqref{FiFj} holds.\\
\textbf{Case~1.}\quad $|i-j|>1$, $B_{ij}=0$.
Suppose $i<j$. From \eqref{kikj}-\eqref{XiXj}, we have
\begin{align*}
&H_{i}^{+}(u)H_{j}^{-}(v)=H_{j}^{-}(v)H_{i}^{+}(u),\\
&H_{i}^{-}(u)H_{i}^{+}(v)=H_{i}^{+}(v)H_{i}^{-}(u),\\
&H_{i}^{\pm}(u)^{-1}E_{j}(v)H_{i}^{\pm}(u)=E_{j}(v),\\
&H_{i}^{\pm}(u)F_{j}(v)H_{i}^{\pm}(u)^{-1}=F_{j}(v),\\
&E_{i}(u)E_{j}(v)=(-1)^{|\alpha_{i}||\alpha_{j}|}E_{j}(v)E_{i}(u),\\
&F_{i}(u)F_{j}(v)=(-1)^{|\alpha_{i}||\alpha_{j}|}F_{j}(v)F_{i}(u).
\end{align*}
\textbf{Case~2.}\quad $i=j$, $B_{ij}=\frac{1}{2}(d_{i}+d_{i+1})$.
From \eqref{kikj}-\eqref{XiXj}, we get
\begin{align*}
&H_{i}^{+}(u)H_{i}^{-}(v)=\frac{\left(u_{-}-v_{+}-d_{i+1}h\right)\left(u_{+}-v_{-}
+d_{i}h\right)}{\left(u_{-}-v_{+}+d_{i}h\right)\left(u_{+}-v_{-}-d_{i+1}h\right)}H_{i}^{-}(v)H_{i}^{+}(u),\\
&H_{i}^{-}(u)H_{i}^{+}(v)=\frac{\left(u_{-}-v_{+}+d_{i+1}h\right)\left(u_{+}-v_{-}-d_{i}h\right)}{\left(u_{+}-v_{-}+d_{i+1}h\right)\left(u_{-}-v_{+}-d_{i}h\right)}
H_{i}^{+}(v)H_{i}^{-}(u),\\
&H_{i}^{\pm}(u)^{-1}E_{i}(v)H_{i}^{\pm}(u)=\frac{u_{\pm}-v-d_{i+1}h}{u_{\pm}-v+d_{i}h}E_{i}(v),\\
&H_{i}^{\pm}(u)F_{i}(v)H_{i}^{\pm}(u)^{-1}
=\frac{u_{\mp}-v-d_{i+1}h}{u_{\mp}-v+d_{i}h}F_{i}(v),\\
&E_{i}(u)E_{i}(v)=(-1)^{|\alpha_{i}|}\frac{u-v+d_{i+1}h}{u-v-d_{i}h}E_{i}(v)E_{i}(u),\\
&F_{i}(u)F_{i}(v)=(-1)^{|\alpha_{i}|}\frac{u-v-d_{i+1}h}{u-v+d_{i}h}F_{i}(v)F_{i}(u).
\end{align*}
\textbf{Case~3.}\quad $j=i+1$, $B_{ij}=-\frac{1}{2}d_{i+1}$.
From \eqref{kikj}-\eqref{XiXj}, we can deduce
\begin{align*}
&H_{i}^{+}(u)H_{i+1}^{-}(v)=\frac{\left(u_{-}-v_{+}+\frac{1}{2}d_{i+1}h\right)\left(u_{+}-v_{-}-\frac{1}{2}d_{i+1}h\right)}{\left(u_{-}-v_{+}-\frac{1}{2}d_{i+1}h\right)\left(u_{+}-v_{-}+
\frac{1}{2}d_{i+1}h\right)}H_{i+1}^{-}(v)H_{i}^{+}(u),\\
&H_{i}^{-}(u)H_{i+1}^{+}(v)=\frac{\left(u_{+}-v_{-}+\frac{1}{2}d_{i+1}h\right)\left(u_{-}-v_{+}-\frac{1}{2}
d_{i+1}h\right)}{\left(u_{+}-v_{-}-\frac{1}{2}d_{i+1}h\right)\left(u_{-}-
v_{+}+\frac{1}{2}d_{i+1}h\right)}H_{i+1}^{+}(v)H_{i}^{-}(u),\\
&H_{i}^{\pm}(u)^{-1}E_{i+1}(v)H_{i}^{\pm}(u)=\frac{u_{\pm}-v+\frac{1}{2}d_{i+1}h}{u_{\pm}-v-\frac{1}{2}d_{i+1}h}E_{i+1}(v),\\
&H_{i}^{\pm}(u)F_{i+1}(v)H_{i}^{\pm}(u)^{-1}=\frac{u_{\mp}-v+\frac{1}{2}d_{i+1}h}{u_{\mp}-v-\frac{1}{2}d_{i+1}h}F_{i+1}(v),\\
&E_{i}(u)E_{i+1}(v)=(-1)^{|\alpha_{i}||\alpha_{i+1}|}\frac{u-v-\frac{1}{2}d_{i+1}h}{u-v+\frac{1}{2}
d_{i+1}h}E_{i+1}(v)E_{i}(u),\\
&F_{i}(u)F_{i+1}(v)=(-1)^{|\alpha_{i}||\alpha_{i+1}|}\frac{u-v+\frac{1}{2}d_{i+1}h}{u-v-\frac{1}{2}
d_{i+1}h}F_{i+1}(v)F_{i}(u).
\end{align*}
\textbf{Case~4.}\quad~$j=i-1$, $B_{ij}=-\frac{1}{2}d_{i}$.
From \eqref{kikj}-\eqref{XiXj}, we derive
\begin{align*}
&H_{i}^{-}(u)H_{i-1}^{+}(v)=\frac{\left(u_{+}-v_{-}+\frac{1}{2}d_{i}h\right)\left(u_{-}-v_{+}-\frac{1}{2}d_{i}h\right)}{\left(u_{+}-v_{-}-\frac{1}{2}d_{i}h\right)\left(u_{-}-v_{+}+
\frac{1}{2}d_{i}h\right)}H_{i-1}^{+}(v)H_{i}^{-}(u),\\
&H_{i}^{+}(u)H_{i-1}^{-}(v)=\frac{\left(u_{+}-v_{-}-\frac{1}{2}d_{i}h\right)\left(u_{-}-v_{+}+\frac{1}{2}d_{i}h\right)}{\left(u_{+}-v_{-}+\frac{1}{2}d_{i}h\right)\left(u_{-}-v_{+}-\frac{1}{2}d_{i}h\right)}H_{i-1}^{-}(v)H_{i}^{+}(u),\\
&H_{i}^{\pm}(u)^{-1}E_{i-1}(v)H_{i}^{\pm}(u)=\frac{u_{\pm}-v+\frac{1}{2}d_{i}h}{u_{\pm}-v-\frac{1}{2}d_{i}h}E_{i-1}(v),\\
&H_{i}^{\pm}(u)F_{i-1}(v)H_{i}^{\pm}(u)^{-1}=\frac{u_{\mp}-v+\frac{1}{2}d_{i}h}{u_{\mp}-v-\frac{1}{2}d_{i}h}F_{i-1}(v),\\
&E_{i}(u)E_{i-1}(v)=(-1)^{|\alpha_{i}||\alpha_{i-1}|}\frac{u-v-\frac{1}{2}d_{i}h}{u-v+\frac{1}{2}
d_{i}h}E_{i-1}(v)E_{i}(u),\\
&F_{i}(u)F_{i-1}(v)=(-1)^{|\alpha_{i}||\alpha_{i-1}|}\frac{u-v+\frac{1}{2}d_{i}h}{u-v-\frac{1}{2}
d_{i}h}F_{i-1}(v)F_{i}(u).
\end{align*}

Combining the above four cases, we derive the relation \eqref{Hi+Hj-}-\eqref{FiFj}.

From \eqref{Xi+Xj-}, we get
\begin{align*}
\left[E_{i}(u),F_{j}(v)\right]=\frac{\delta_{ij}}{h}\{\delta(\frac{v_{+}}{u_{-}})H_{i}^{+}(u_{-})-\delta(\frac{v_{-}}{u_{+}})H_{i}^{-}(v_{-})\}.
\end{align*}

From \eqref{serre  i-1} and \eqref{serre i-2}, we get
\begin{align*}
&{\rm Sym}_{u_{1},u_{2}}\left[E_{i}(u_{1}),\left[E_{i}(u_{2}),E_{i\pm 1}(v)\right]\right]=0,\quad |\alpha_{i}|=0,\\
&{\rm Sym}_{u_{1},u_{2}}\left[F_{i}(u_{1}),\left[F_{i}(u_{2}),F_{i\pm 1}(v)\right]\right]=0,\quad |\alpha_{i}|=0,\\
&{\rm Sym}_{u_{1},u_{2}}\left[E_{i}(u_{1}),\left[E_{i+1}(v_{1}),\left[E_{i}(u_{2}),
E_{i-1}(v_{2})\right]\right]\right]=0,\quad |\alpha_{i}|=1,\\
&{\rm Sym}_{u_{1},u_{2}}\left[F_{i}(u_{1}),\left[F_{i+1}(v_{1}),\left[F_{i}(u_{2}),
F_{i-1}(v_{2})\right]\right]\right]=0,\quad |\alpha_{i}|=1.
\end{align*}

This completes the proof.
\end{proof}

Define $\widehat{\mathrm{DY}}_{h}(\mathfrak{sl}^{\mathfrak{s}}_{m|n})$ to be the algebra generated by the coefficients of the series $H_{i}^{\pm}(u)$, $E_{i}(u)$, and $F_{i}(u)$ for $i\in I_{\mathfrak{s}}'$, together with a central element $C$, subject to the relations stated in Lemma \ref{Drinfeld generator sl}. We shall prove that $\widehat{\mathrm{DY}}_{h}(\mathfrak{sl}^{\mathfrak{s}}_{m|n})$ is isomorphic to $\mathrm{DY}_{h}(\mathfrak{sl}^{\mathfrak{s}}_{m|n})$; this isomorphic presentation will be called the Drinfeld presentation of the double Yangian for $\mathfrak{sl}^{\mathfrak{s}}_{m|n}$.

\begin{theo}\label{Iso-sl}
The associative superalgebras $\mathrm{DY}_{h}(\mathfrak{sl}^{\mathfrak{s}}_{m|n})$ and $\widehat{\mathrm{DY}}_{h}(\mathfrak{sl}^{\mathfrak{s}}_{m|n})$ are isomorphic.
\end{theo}
\begin{proof}
Lemma \ref{SL RTT} and lemma \ref{Drinfeld generator sl} imply that there is a surjective homomorphism $\varphi:\widehat{\mathrm{DY}}_{h}(\mathfrak{sl}^{\mathfrak{s}}_{m|n})\rightarrow \mathrm{DY}_{h}(\mathfrak{sl}^{\mathfrak{s}}_{m|n})$, which takes each generator in $\widehat{\mathrm{DY}}_{h}(\mathfrak{sl}^{\mathfrak{s}}_{m|n})$ to the element of the same name in the $\mathrm{DY}_{h}(\mathfrak{sl}^{\mathfrak{s}}_{m|n})$. The same arguments as in the proof of Theorem \ref{iso gl} show that $\varphi$ is also injective. This completes the proof.
\end{proof}

\medskip
\section{Bosonization of level-1 modules}

In this section, the level-1 modules for $\mathrm{DY}_{h}(\mathfrak{g})$ are constructed in terms of bosons for $\mathfrak{g}=\mathfrak{gl}^{\mathfrak{s}}_{m|n}$ and $\mathfrak{sl}^{\mathfrak{s}}_{m|n}$, where $\mathfrak{s}\in\mathcal{S}_{m|n}$. The construction is based on the Drinfeld presentation and follows the ideas of \cite{I96} and \cite{JYL20}. We introduce bosonic oscillators $\{a_{i,k} \mid i\in I_{\mathfrak{s}}, ~k\in\mathbb{Z}^{*}\}$ satisfying
$$
[a_{i,k},a_{j,l}]=d_{i}k\delta_{ij}\delta_{k+l,0}.
$$

\begin{rema}
Throughout this section, the bracket $[\cdot,\cdot]$ refers to the bracket operation for ordinary (non-super) Lie algebras.
\end{rema}

\subsection{The $\mathfrak{gl}^{\mathfrak{s}}_{m|n}$ case}

Let $\bar{\mathcal{Q}}_{\mathfrak{s}}=\mathcal{Q}_{\mathfrak{s}}+\sum_{i\in I_{\mathfrak{s}}}\mathbb{Z}\epsilon_{i}$ and $\mathcal{A}[\bar{\mathcal{Q}}_{\mathfrak{s}}]$ denote the group algebra of the root lattice $\bar{\mathcal{Q}}_{\mathfrak{s}}$ over $\mathcal{A}$, with the basis $\{e^{\alpha}~|~\alpha\in \bar{\mathcal{Q}}_{\mathfrak{s}}\}$. We introduce the Fock space $\mathcal{F}_{\mathfrak{s}}=\mathcal{A}[a_{i,-k}~(i \in I_{\mathfrak{s}}, k\in\mathbb{Z}_{+})] \widetilde{\otimes}\mathcal{A}[\bar{\mathcal{Q}}_{\mathfrak{s}}]$. On this space, we define the action of the operators $a_{i,k}, \partial \epsilon_{i}, e^{\epsilon_{i}}~(i \in I_{\mathfrak{s}})$ by
\begin{align*}
a_{i,k}\cdot f\otimes e^{\beta}&=\left\{\begin{array}{ll}
a_{i,k}f\otimes e^{\beta}, & \text{if } k<0, \\
 \left[a_{i,k},f\right] \otimes e^{\beta}, & \text{if } k>0;
\end{array}
 \right.\\
\partial \epsilon_{i}\cdot f\otimes e^{\beta}&=\epsilon_{i}(h_{\beta}) f\otimes e^{\beta},\\
e^{\epsilon_{i}}\cdot f\otimes e^{\beta}&=c(\epsilon_{i},\beta)f\otimes e^{\epsilon_{i}+\beta},
\end{align*}
for $f\otimes e^{\beta}\in \mathcal{F}_{\mathfrak{s}}$. Here, $c(\cdot,\cdot)$ is a 2-cocycle on $\bar{\mathcal{Q}}_{\mathfrak{s}}$ satisfying
\begin{align*}
c(\alpha,\beta)&=(-1)^{|\alpha||\beta|}c(\beta,\alpha),\\
c(\alpha,\beta)c(\alpha+\beta,\gamma)&=c(\beta,\gamma)c(\alpha,\beta+\gamma),
\end{align*}
for $\alpha,\beta,\gamma\in \bar{\mathcal{Q}}_{\mathfrak{s}}$.

\begin{theo}
For $i \in I_{\mathfrak{s}}, j \in I_{\mathfrak{s}}'$, the following assignment defines a $\mathrm{DY}_{h}(\mathfrak{gl}^{\mathfrak{s}}_{m|n})$-module structure on $\mathcal{F}_{\mathfrak{s}}$.
{\allowdisplaybreaks
\begin{align*}
k_{i}^{+}(u)\mapsto \exp &\left(-\sum_{k>0}\frac{a_{i,k}}{k}\left(\left(u+
\frac{1}{2}d_{i}h\right)^{-k}-\left(u-\frac{1}{2}d_{i}h\right)^{-k}\right)\right)
\left(\frac{u-\frac{1}{2}h}{u+\frac{1}{2}h}\right)^{\partial \epsilon_{i}},\\
k_{i}^{-}(u)\mapsto \exp &\left(\sum_{k>0,r>i}\frac{a_{i,-k}}{k}\left(\left(u+d_{i}h\right)^{k}+\left(u-d_{i}h\right)^{k}-2u^{k}\right)+\sum_{k>0}\frac{a_{i,-k}}{k}\left(\left(u-d_{i}h\right)^{k}-u^{k}\right)\right),\\
\frac{1}{h}X_{j}^{+}(u) \mapsto \exp &\left(-\sum_{k>0}d_{j} \frac{a_{j,-k}}{k}\left(u-\frac{1}{2}d_{j}h-\frac{1}{4}h\right)^{k}+\sum_{k>0}d_{j+1}\frac{a_{j+1,-k}}{k}\left(u+\frac{1}{2}d_{j+1}h
-\frac{1}{4}h\right)^{k}\right)\\
\exp &\left(\sum_{k>0}\frac{d_{j}a_{j,k}-d_{j+1}a_{j+1,k}}{k}
\left(u+\frac{1}{4}h\right)^{-k}\right)e^{\alpha_{j}}\left((-1)^{j-1}
\left(u+\frac{1}{4}h\right)\right)^{\partial \alpha_{j}},\\
\frac{1}{h}X_{j}^{-}(u)\mapsto \exp &\left(\sum_{k>0}d_{j}
\frac{a_{j,-k}}{k}\left(u-\frac{1}{2}d_{j}h+\frac{1}{4}h\right)^{k}-\sum_{k>0}d_{j+1}\frac{a_{j+1,-k}}{k}\left(u+\frac{1}{2}d_{j+1}h
+\frac{1}{4}h\right)^{k}\right)\\
\exp &\left(\sum_{k>0}\frac{-d_{j}a_{j,k}+d_{j+1}a_{j+1,k}}{k}
\left(u-\frac{1}{4}h\right)^{-k}\right)e^{-\alpha_{j}}\left((-1)^{j-1}
\left(u-\frac{1}{4}h\right)\right)^{-\partial \alpha_{j}}.
\end{align*}
}
\end{theo}
\begin{proof}
For $\alpha,\beta,\gamma \in \bar{\mathcal{Q}}_{\mathfrak{s}}$, the definition of $\partial \alpha, e^{\beta}$ and $e^{\gamma}$ implies that $[\partial \alpha, e^{\beta}]=\alpha(h_{\beta})e^{\beta}$ and $e^{\beta}e^{\gamma}=(-1)^{|\beta||\gamma|}e^{\gamma}e^{\beta}$. Then we have
\begin{align*}
\left(\frac{u-\frac{1}{2}h}{u+\frac{1}{2}h}\right)^{\partial \epsilon_{i}}
e^{\alpha_{j}}&=e^{\alpha_{j}}\left(\frac{u-\frac{1}{2}h}{u+\frac{1}{2}h}\right)^{\partial \epsilon_{i}}\left(\frac{u-\frac{1}{2}h}{u+\frac{1}{2}h}\right)^{\epsilon_{i}(h_{\alpha_{j}})},\\
e^{\alpha_{i}}\left((-1)^{j-1}
\left(u+\frac{1}{4}h\right)\right)^{-\partial \alpha_{j}}&=\left((-1)^{j-1}
\left(u+\frac{1}{4}h\right)\right)^{-\partial \alpha_{j}}e^{\alpha_{i}}\left((-1)^{j-1}
\left(u+\frac{1}{4}h\right)\right)^{\alpha_{i}(h_{\alpha_{j}})},\\
\left((-1)^{j-1}
\left(u-\frac{1}{4}h\right)\right)^{\partial \alpha_{i}}e^{-\alpha_{j}}&=e^{-\alpha_{j}}\left((-1)^{j-1}
\left(u-\frac{1}{4}h\right)\right)^{\partial \alpha_{i}}\left((-1)^{j-1}
\left(u-\frac{1}{4}h\right)\right)^{-\alpha_{i}(h_{\alpha_{j}})}.
\end{align*}

Next, we need to verify the relations in Theorem \ref{Drinfeld relation 1} hold when $C=1$. We now verify the relation \eqref{XiXi}. Recall that if $[A,B]$ commutes with $A$ and $B$, then $\exp (A)\exp (B)=\exp (B)\exp (A)\exp ([A,B])$. Using the commutation relations $[a_{i,k},a_{j,l}]=d_{i}k\delta_{i,j}\delta_{k+l,0}$ and the expansion $\ln (1-u)=-\sum_{k>0}\frac{u^{k}}{k}$, we derive the following relations:
\begin{align*}
&\exp \left(-\sum_{k>0}\frac{d_{j}a_{j,-k}}{k}\left(u-\frac{1}{2}d_{j}h-\frac{1}{4}h\right)^{k}\right)\exp \left(\sum_{k>0}\frac{d_{j}a_{j,k}}{k}\left(v+\frac{1}{4}h\right)^{-k}\right)\\
=~&\exp \left(\sum_{k>0}\frac{d_{j}a_{j,k}}{k}\left(v+\frac{1}{4}h\right)^{-k}\right)\exp \left(-\sum_{k>0}\frac{d_{j}a_{j,-k}}{k}\left(u-\frac{1}{2}d_{j}h-\frac{1}{4}h\right)^{k}\right)\left(\frac{v+\frac{1}{4}h}{v-u+\frac{1+d_{j}}{2}h}\right)^{d_{j}},\\
&\exp \left(\sum_{k>0}\frac{d_{j+1}a_{j+1,-k}}{k}\left(u+\frac{1}{2}d_{j+1}h-\frac{1}{4}h\right)^{k}\right)\exp \left(-\sum_{k>0}\frac{d_{j+1}a_{j+1,k}}{k}\left(v+\frac{1}{4}h\right)^{-k}\right)\\
=~&\exp \left(-\sum_{k>0}\frac{d_{j+1}a_{j+1,k}}{k}\left(v+\frac{1}{4}h\right)^{-k}\right)\exp \left(\sum_{k>0}\frac{d_{j+1}a_{j+1,-k}}{k}\left(u+\frac{1}{2}d_{j+1}h-\frac{1}{4}h\right)^{k}\right)\\
&\cdot \left(\frac{v+\frac{1}{4}h}{v-u+\frac{1-d_{j+1}}{2}h}\right)^{d_{j+1}},\\
&\exp \left(\sum_{k>0}\frac{d_{j}a_{j,k}}{k}\left(u+\frac{1}{4}h\right)^{-k}\right)\exp \left(-\sum_{k>0}\frac{d_{j}a_{j,-k}}{k}\left(v-\frac{1}{2}d_{j}h-\frac{1}{4}h\right)^{k}\right)\\
=~&\exp \left(-\sum_{k>0}\frac{d_{j}a_{j,-k}}{k}\left(v-\frac{1}{2}d_{j}h-\frac{1}{4}h\right)^{k}\right)\exp \left(\sum_{k>0}\frac{d_{j}a_{j,k}}{k}\left(u+\frac{1}{4}h\right)^{-k}\right)\left(\frac{u-v+\frac{1+d_{j}}{2}h}{u+\frac{1}{4}h}\right)^{d_{j}},\\
&\exp \left(-\sum_{k>0}\frac{d_{j+1}a_{j+1,k}}{k}\left(u+\frac{1}{4}h\right)^{-k}\right)\exp \left(\sum_{k>0}\frac{d_{j+1}a_{j+1,-k}}{k}\left(v+\frac{1}{2}d_{j+1}h-\frac{1}{4}h\right)^{k}\right)\\
=~&\exp \left(\sum_{k>0}\frac{d_{j+1}a_{j+1,-k}}{k}\left(v+\frac{1}{2}d_{j+1}h-\frac{1}{4}h\right)^{k}\right)\exp \left(-\sum_{k>0}\frac{d_{j+1}a_{j+1,k}}{k}\left(u+\frac{1}{4}h\right)^{-k}\right)\\
&\cdot \left(\frac{u-v+\frac{1-d_{j+1}}{2}h}{u+\frac{1}{4}h}\right)^{d_{j+1}}.
\end{align*}

Since $e^{\alpha_{j}}e^{\alpha_{j}}=(-1)^{|\alpha_{j}|}e^{\alpha_{j}}e^{\alpha_{j}}$ and 
\begin{align*}
&\left((-1)^{j-1}
\left(u+\frac{1}{4}h\right)\right)^{\alpha_{i}(h_{\alpha_{j}})}\left((-1)^{j-1}
\left(v+\frac{1}{4}h\right)\right)^{-\alpha_{i}(h_{\alpha_{j}})}\left(\frac{v+\frac{1}{4}h}{v-u+\frac{1+d_{j}}{2}h}\right)^{d_{j}}\\
&\left(\frac{v+\frac{1}{4}h}{v-u+\frac{1-d_{j+1}}{2}h}\right)^{d_{j+1}}\left(\frac{u-v+\frac{1+d_{j}}{2}h}{u+\frac{1}{4}h}\right)^{d_{j}}\left(\frac{u-v+\frac{1-d_{j+1}}{2}h}{u+\frac{1}{4}h}\right)^{d_{j+1}}\\
=&\left(-1\right)^{d_{j}+d_{j+1}}\left(\frac{u-v+\frac{1+d_{j}}{2}h}{u-v-\frac{1+d_{j}}{2}h}\right)^{d_{j}}\left(\frac{u-v+\frac{1-d_{j+1}}{2}h}{u-v-\frac{1-d_{j+1}}{2}h}\right)^{d_{j+1}},
\end{align*}
it follows that
\begin{align*}
  X_{j}^{+}(u)X_{j}^{+}(v)=  \left\{
  \begin{array}{ll}
    \frac{u-v+h}{u-v-h}X_{j}^{+}(v)X_{j}^{+}(u), & d_{j}=d_{j+1}=1, \\
    \frac{u-v-h}{u-v+h}X_{j}^{+}(v)X_{j}^{+}(u), & d_{j}=d_{j+1}=-1,\\
    -X_{j}^{+}(v)X_{j}^{+}(u), & d_{j}+d_{j+1}=0.
  \end{array}
\right.
\end{align*}
We therefore obtain the operator product expansion:
\begin{align*}
\left\{
  \begin{array}{ll}
    X_{j}^{+}(u)X_{j}^{+}(v)=\frac{u-v+d_{j+1}h}{u-v-d_{j}h}X_{j}^{+}(v)X_{j}^{+}(u), & |\alpha_{j}|=0, \\
    X_{j}^{+}(u)X_{j}^{+}(v)=-X_{j}^{+}(v)X_{j}^{+}(u), & |\alpha_{j}|=1,
  \end{array}
\right.
\end{align*}
which implies that \eqref{XiXi} holds. The remaining relations in Theorem \ref{Drinfeld relation 1} follow by similar calculations.
\end{proof}

\subsection{The $\mathfrak{sl}^{\mathfrak{s}}_{m|n}$ case}

We keep the notations from the $\mathfrak{gl}^{\mathfrak{s}}_{m|n}$ case unless otherwise stated. Set
$$
\mathcal{F}_{\mathfrak{s}}:=\mathcal{A}[a_{i,-k}(i \in I_{\mathfrak{s}}', k\in\mathbb{Z}_{+})]\widetilde{\otimes}\mathcal{A}[\mathcal{Q}_{\mathfrak{s}}],
$$
and let the operators $a_{i,k}, \partial \alpha_{i}, e^{\alpha_{i}}~( i \in I_{\mathfrak{s}}')$ act on it as in the previous subsection. We impose the convention $a_{0,-k}=a_{N,-k}=0$.

\begin{theo}\label{Vertex sl}
For $i \in I_{\mathfrak{s}}'$, the following assignment defines a $\mathrm{DY}_{h}(\mathfrak{sl}^{\mathfrak{s}}_{m|n})$-module structure on $\mathcal{F}_{\mathfrak{s}}$.
\begin{align*}
H_{i}^{+}(u)\mapsto \exp &\left(-\sum_{k>0}
\frac{d_{i}a_{i,k}}{k}\left(\left(u+\frac{1}{2}h\right)^{-k}-\left(u-\frac{1}{2}h
\right)^{-k}\right)\right)\left(\frac{u-\frac{1}{2}h}{u+\frac{1}{2}h}\right)^{-\partial \alpha_{i}},\\
H_{i}^{-}(u)\mapsto \exp &\left(-\sum_{k>0}\frac{a_{i,-k}}{k}\left(\left(u+d_{i+1}h\right)^{k}-
\left(u-d_{i}h\right)^{k}\right)\right.\\
+&\left.\sum_{k>0}\frac{a_{i+1,-k}}{k}\left(\left(u+\frac{d_{i+1}}{2}h\right)^{k}-
\left(u-\frac{d_{i+1}}{2}h\right)^{k}\right)\right.\\
+&\left.\sum_{k>0}\frac{a_{i-1,-k}}{k}\left(\left(u+\frac{d_{i}}{2}h\right)^{k}-
\left(u-\frac{d_{i}}{2}h\right)^{k}\right)\right),\\
E_{i}(u)\mapsto \exp &\left(\sum_{k>0}
\frac{d_{i}a_{i,-k}}{k}\left(\left(u+\left(\frac{d_{i+1}}{2}h
-\frac{1}{4}h\right)\right)^{k}+(-1)^{|\alpha_{i}|}\left(u-\left(\frac{d_{i}}{2}h
+\frac{1}{4}h\right)\right)^{k}\right)\right.\\
-&\left.\sum_{k>0}\frac{d_{i}a_{i-1,-k}}{k}\left(u-\frac{1}{4}h\right)^{k}-
\sum_{k>0}\frac{d_{i+1}a_{i+1,-k}}{k}\left(u-\frac{1}{4}h\right)^{k}\right)\\
\exp &\left(-\sum_{k>0}\frac{d_{i}a_{i,k}}{k}
\left(u+\frac{1}{4}h\right)^{-k}\right)e^{ \alpha_{i}}\left((-1)^{i-1}\left(u+\frac{1}{4}h\right)\right)^{\partial \alpha_{i}},\\
F_{i}(u)\mapsto \exp &\left(-\sum_{k>0}\frac{d_{i}a_{i,-k}}{k}
\left(\left(u+\left(\frac{d_{i+1}}{2}h+\frac{1}{4}h\right)\right)^{k}+(-1)^{|\alpha_{i}|}
\left(u-\left(\frac{d_{i}}{2}h-\frac{1}{4}h\right)\right)^{k}\right)\right.\\
+&\left.\sum_{k>0}\frac{d_{i+1}a_{i+1,-k}}{k}\left(u+\frac{1}{4}h\right)^{k}+
\sum_{k>0}\frac{d_{i}a_{i-1,-k}}{k}\left(u+\frac{1}{4}h\right)^{k}\right)\\
\exp &\left(\sum_{k>0}\frac{d_{i}a_{i,k}}{k}
\left(u-\frac{1}{4}h\right)^{-k}\right)e^{-\alpha_{i}}\left((-1)^{i-1}\left(u-
\frac{1}{4}h\right)\right)^{-\partial \alpha_{i}}.
\end{align*}
\end{theo}
\begin{proof}
We need to verify the relations in Lemma \ref{Drinfeld generator sl} with $C=1$. These relations can be checked in the same way as in the $\mathfrak{gl}_{m|n}^{\mathfrak{s}}$ case.
\end{proof}
\begin{rema}
In the case $n=0$, the relations of Theorem \ref{Vertex sl} are identical to those in $\mathrm{DY}_{h}(\mathfrak{sl}_{m})$ presented in Iohara's work \cite{I96}.
\end{rema}

\medskip
\noindent\textbf{Acknowledgments} H. Lin is supported by the Postdoctoral Fellowship Program of CPSF (GZC20252014). H. Zhang is supported by the National Natural Science Foundation of China 12271332.


\begin{thebibliography}{99}
\vskip5pt
\bibitem{AACFR03}
 D. Arnaudon, J. Avan, N. Crampé, L. Frappat, E. Ragoucy, R-matrix presentation for super-Yangians $Y(\mathfrak{osp}(m|2n))$, J. Math. Phys.
\textbf{44} (2003), 302-308.

\bibitem{BK25-1}
L. Bagnoli, S. Ko\v zi\'c, A note on the quantum Berezinian for the double Yangian of the Lie superalgebra $\mathfrak{gl}_{m|n}$, Algebr. Represent. Theory \textbf{28}(1) (2025) 143-155.

\bibitem{BK25-2}
L. Bagnoli, S. Ko\v zi\'c, Double Yangian and reflection algebras of the Lie
superalgebra {$\mathfrak{gl}_{m|n}$}, Commun. Contemp. Math. \textbf{27}(2) (2025), Paper No. 2450007, 25.

\bibitem{BK05}
J. Brundan, A. Kleshchev, Parabolic presentations of the Yangian $Y(\mathfrak{gl}_n)$, Comm. Math. Phys. \textbf{254}(1) (2005), 191-220.

\bibitem{CH23}
H. Chang, H. Hu, A note on the center of the super Yangian {$Y_{M|N}(\mathfrak{ s})$}, J. Algebra \textbf{633} (2023) 648-665.

\bibitem{CP94}
V. Chari, A. Pressley, A guide to quantum groups, Cambridge University Press, Cambridge (1994), xvi+651.

\bibitem{D85}
V.G. Drinfeld, Hopf algebras and the quantum Yang-Baxter equation, Dokl. Akad. Nauk SSSR \textbf{283}(5) (1985) 1060-1064.

\bibitem{D87}
V. G. Drinfeld, Quantum groups, Proceedings of the International Congress of Mathematicians, Vol. 1, 2 (Berkeley, Calif., 1986), Amer. Math. Soc., Providence, RI, (1987), 798-820.

\bibitem{D88}
V.G. Drinfeld, A new realization of Yangians and of quantum affine algebras, Dokl. Akad. Nauk SSSR \textbf{296}(1) (1987) 13-17.

\bibitem{DHHZ98}
 X. M. Ding, B. Y. Hou, B. Yuan Hou, L. Zhao, Free boson representation of $\mathrm{DY}_{h}(\mathfrak{gl}_{N})_{k}$, J. Math. Phys. \textbf{39}(4) (1998), 2273-2289.

\bibitem{EFM01}
E. Frenkel, E. Mukhin, Combinatorics of $q$-characters of finite-dimensional representations of quantum affine algebras, Comm. Math. Phys. \textbf{216}(1) (2001), 23-57.

\bibitem{FJ88}
I. B. Frenkel, N. Jing, Vertex representations of quantum affine algebras, Proc. Nat. Acad. Sci. U.S.A. \textbf{85}(24) (1988), 9373-9377.

\bibitem{BFR92}
I. B. Frenkel, N. Y. Reshetikhin, Quantum affine algebras and holonomic difference equations, Comm. Math. Phys. \textbf{146}(1) (1992), 1-60.

\bibitem{EFR96}
E. Frenkel, N. Reshetikhin, Quantum affine algebras and deformations of the Virasoro and $\mathcal{W}$-algebras, Comm. Math. Phys. \textbf{178}(1) (1996), 237-264.

\bibitem{FRT90}
L. Faddeev, N. Reshetikhin, L. Takhtajan, Quantization of Lie groups and Lie algebras, Leningrad Math. J. \textbf{1}(1) (1990), 193–225.

\bibitem{GGRW05}
I. Gelfand, S. Gelfand, V. Retakh, R.L. Wilson, Quasideterminants, Adv. Math. \textbf{193}(1) (2005), 56-141.

\bibitem{G07}
L. Gow, Gauss decomposition of the Yangian {$Y({\mathfrak{gl}}_{m|n})$}, Comm. Math. Phys. \textbf{276}(3) (2007), 799-825.

\bibitem{I96}
K. Iohara, Bosonic representations of Yangian double $\mathcal{D}Y_{\hbar}(\mathfrak{g})$ with $\mathfrak{g}=\mathfrak{gl}_{N}$, $\mathfrak{sl}_{N}$, J. Phys. A \textbf{29}(15) (1996), 4593-4621.

\bibitem{Ji85}
M. Jimbo, A $q$-difference analogue of $U(g)$ and the Yang-Baxter equation, Lett. Math. Phys. \textbf{10}(1) (1985) 63-69.

\bibitem{JLM18}
N. Jing, M. Liu, A. Molev, Isomorphism between the $R$-matrix and Drinfeld presentations of Yangian in types~$B$, $C$ and~$D$, Comm. Math.
 Phys. \textbf{361}(3) (2018), 827-872.

\bibitem{JYL20}
N. Jing, F. Yang, M. Liu, Yangian doubles of classical types and their vertex representations, J. Math. Phys. \textbf{61}(5) (2020), 051704, 39.

\bibitem{K77}
V.G. Kac, Lie superalgebras, Adv. Math. \textbf{26} (1977), 8-96.

\bibitem{K97}
S. Khoroshkin, Central extension of the Yangian double, Alg$\acute{e}$bre non commutative, groupes quantiques et invariants (Reims, 1995), S$\acute{e}$min. Congr., vol. 2, Soc. Math. France, Paris, (1997), 119-135.

\bibitem{Ko18}
S. Ko$\check{z}$i$\acute{c}$, Commutative operators for double Yangian ${\rm DY}(\mathfrak{sl}_n)$, Glas. Mat. Ser. III \textbf{53}(73) (2018), 97-113.

\bibitem{L22}
K. Lu, A note on odd reflections of super Yangian and Bethe ansatz, Lett. Math. Phys. \textbf{112}(2) (2022), Paper No. 29, 26.

\bibitem{LSS}
D. Leites, M. Saveliev, V. Serganova, Embeddings of $\mathfrak{osp}(N/2)$ and the associated nonlinear supersymmetric equations, Group theoretical methods in physics, Vol. I (Yurmala, 1985), VNU Sci. Press, Utrecht (1986) 255-297.

\bibitem{LZ25}
H. Lin, H. Zhang, Representations of Quantum Affine General Linear Superalgebras at Arbitrary 01-Sequences, arXiv preprint, arXiv:2511.02393.

\bibitem{M22}
A. Molev, Odd reflections in the {Y}angian associated with $\mathfrak{gl}(m|n)$, Lett. Math. Phys. \textbf{112}(1) (2022), Paper No. 8, 15.

\bibitem{M24}
A. Molev, A Drinfeld-type presentation of the orthosymplectic Yangians, Algebr. Represent. Theory \textbf{27}(1) (2024), 469-494.

\bibitem{MNO96}
A. Molev, M. Nazarov, G. Ol'shanskii, Yangians and classical Lie algebras, Russ. Math. Surv. \textbf{51}(2) (1996), 205–282.

\bibitem{Mu12}
I.M. Musson, Lie superalgebras and enveloping algebras, Grad. Stud. Math., \textbf{131} American Mathematical Society, Providence, RI, (2012), xx+488 pp.

\bibitem{N91}
M. Nazarov, Quantum Berezinian and the classical Capelli identity, Lett. Math. Phys. \textbf{21}(2) (1991), 123-131.

\bibitem{P16}
Y. Peng, Parabolic presentations of the super Yangian $Y(\mathfrak{gl}_{M|N})$ associated with arbitrary 01-sequences, Comm. Math. Phys. \textbf{346}(1) (2016) 313-347.

\bibitem{R90}
N. Y. Reshetikhin, Quasitriangular Hopf algebras and invariants of links, Leningrad Math. J. \textbf{1}(2) (1990), 491–513.

\bibitem{RT91}
N. Reshetikhin, V. G. Turaev, Invariants of $3$-manifolds via link polynomials and quantum groups, Invent. Math. \textbf{103}(3) (1991), 547-597.

\bibitem{Ta84}
V. O. Tarasov, The structure of quantum $L$-operators for the $R$-matrix of the $XXZ$-model, Teoret. Mat. Fiz. \textbf{61}(2) (1984), 163-173.

\bibitem{T20}
A. Tsymbaliuk, Shuffle algebra realizations of type $A$ super Yangians and quantum affine superalgebras for all Cartan data, Lett. Math. Phys. \textbf{110}(8) (2020), 2083-2111.

\bibitem{XLZ}
P. Xu, H. Lin, H. Zhang, Isomorphism between the $R$-matrix and Drinfeld presentations of quantum affine superalgebra for type $\boldsymbol{A}$, arXiv:2609.35067.

\bibitem{YJ24}
F. Yang, N. Jing, Center of the Yangian double in type $A$, Sci. China Math. \textbf{67}(9) (2024), 1957-1988.

\bibitem{Y99}
H. Yamane, On defining relations of affine Lie superalgebras and affine quantized universal enveloping superalgebras, Publ. Res. Inst. Math.
 Sci. \textbf{35}(3) (1999), 321-390.

\bibitem{XZ}
Y. Xu, R. B. Zhang, Drinfeld realizations and vertex operator representations of quantum affine superalgebras, arXiv:1802.09702.

\bibitem{Z97}
Y. Zhang, Super-Yangian double and its central extension, Phys. Lett. A \textbf{234}(1) (1997), 20-26.
\end{thebibliography}
\end{document}